\documentclass[a4paper,fleqn]{article}
\usepackage{arxiv}

\makeatletter
\def\ps@pprintTitle{%
	\let\@oddhead\@empty
	\let\@evenhead\@empty
	\def\@oddfoot{}%
	\let\@evenfoot\@oddfoot}
\makeatother

\usepackage{amsmath,amsfonts,mathrsfs}

\usepackage{bm}

\usepackage[title]{appendix}
\usepackage{amsthm}
\usepackage[vlined,boxruled]{algorithm2e}
\usepackage{enumitem}
\usepackage{blindtext}
\usepackage{mdframed}
\usepackage[sort,square,numbers]{natbib}
\usepackage{multicol}
\usepackage[bookmarks=true]{hyperref}
\usepackage[numbered]{bookmark}
\usepackage{cleveref}

\usepackage{dsfont}

\usepackage[justification=centering]{caption}
\usepackage{subcaption}

\usepackage{booktabs}

\usepackage{tikz}
\usetikzlibrary{arrows.meta}
\usetikzlibrary{positioning}
\usetikzlibrary{shapes, decorations.pathreplacing}

\usepackage{array}
\usepackage{multirow,multicol}

\usepackage{float}
\usepackage{pgffor}
\usepackage{caption}
\usepackage{subcaption}

\usepackage[section]{placeins}

\usepackage{chemarrow}

\usepackage{autonum}
\AtBeginDocument{
	\newenvironment{equation*}{\begin{equation}}{\end{equation}}
	\newenvironment{align*}{\begin{align}}{\end{align}}
}

\newtheorem{theorem}{Theorem}[section]
\newtheorem{lemma}[theorem]{Lemma}
\newtheorem{proposition}[theorem]{Proposition}

\newtheorem{assum}[theorem]{Assumption}

\theoremstyle{remark}
\newtheorem{remark}[theorem]{Remark}

\theoremstyle{definition}

\newtheorem{Ex}{Example}

\Crefname{Ex}{Example}{Examples}
\Crefname{Def}{Definition}{Definitions}
\Crefname{Thm}{Theorem}{Theorems}
\Crefname{Rem}{Remark}{Remarks}
\Crefname{Lem}{Lemma}{Lemmas}

\newcounter{defcounter}
\newcommand{\Rone}{\mathbb{R}}
\newcommand{\R}{\mathbb{R}}

\newcommand{\q}[1]{``#1''}

\newcommand{\di}{\,\mathrm{d}}

\newcommand\norm[1]{\left\lVert#1\right\rVert}

\newcommand{\Mk}{\mathrm{MAK}}
\newcommand{\Gd}{\mathrm{Gdnv}}
\newcommand{\Cp}{\mathrm{Cap}}

\newcommand{\RR}{\mathcal{K}}
\newcommand{\rrarrow}[1]{\mathop{\longrightarrow}^{#1}}

\definecolor{myred}{HTML}{B83232}
\definecolor{myblue}{HTML}{28568A}
\definecolor{mygreen}{HTML}{288A42}

\newcommand{\changednew}[1]{{\leavevmode \color{black}{#1}}}
\begin{document}

\raggedbottom

% flatex input: [title.tex]
\let\WriteBookmarks\relax
\def\floatpagepagefraction{1}
\def\textpagefraction{.001}
%\shorttitle{Traffic Reaction Model}
%\shortauthors{Pereira et~al.}
%\begin{frontmatter}

\renewcommand{\shorttitle}{\textsc{Traffic Reaction Model}}

\title{The Traffic Reaction Model:
	A kinetic compartmental approach to road traffic modeling}

\author{{\Large M. Pereira$^{1,2}$\thanks{Corresponding author \\ \textit{E-mail addresses}: \texttt{mike.pereira@minesparis.psl.eu} (M. Pereira), \texttt{kulcsar@chalmers.se} (B. Kulcsár), \texttt{lipgyorgy@scl.sztaki.hu} (G. Lipták),
			\texttt{kovacs.mihaly@itk.ppke.hu} (M. Kovács), \texttt{szederkenyi@itk.ppke.hu} (G.~Szederkényi).\\ \textbf{Acknowledgement}: The  authors  acknowledge  the  helpful  discussions  with  M. Vághy, and the comments of the various reviewers whose expertise helped reshape and develop this work. Thank you for being our disciplined  colleague,  Pinar, some memories will never fade.    B.  Kulcsár  acknowledges the  contribution of Transport Area of Advance  at Chalmers  University  of Technology. The  project  has  been  partially  supported by European Comision(F-ENUAC-2022-0003) and Energimyndigheten (2023-00021) through the project E-LaaS:Energy optimal urban Logistics As A Service.  M. Kovács acknowledges the support of the grant NKFIH TKP2021-NVA-02. G. Szederkényi acknowledges the support of the projects NKFIH 145934 and RRF-2.3.1-21-2022-00006.}, B.~Kulcs\'ar$^{2}$,  Gy. Lipták$^{3,4}$, M. Kov\'acs$^{4}$, G. Szederk\'enyi$^{3,4}$} \\ \\	\vspace{1ex}
	$^{1}$Department of Geosciences and Geoengineering, Mines Paris -- PSL University, Fontainebleau, France.\\ 	\vspace{1ex}
	$^{2}$Department of Electrical Engineering, Chalmers University of Technology, G\"oteborg, Sweden.\\ 	\vspace{1ex}
	$^{3}$Institute for Computer Science and Control (SZTAKI),Budapest, Hungary.\\ 	\vspace{1ex}
	$^{4}$Faculty of Information Technology and Bionics, Pázmány Péter Catholic University, Budapest, Hungary.	
	%% examples of more authors
	%%\AND
	%Coauthor \\
	%Affiliation \\
	%% Address \\
	%% \texttt{email} \\
	%% \And
	%% Coauthor \\
	%% Affiliation \\
	%% Address \\
	%% \texttt{email} \\
	%% \And
	%% Coauthor \\
	%% Affiliation \\
	%% Address \\
	%% \texttt{email} \\
}

%\author[2]{M. Pereira}
%\cormark[1]
%\ead{mike.pereira@chalmers.se}
%
%\author[3]{B.~Kulcs\'ar}
%\ead{kulcsar@chalmers.se}
%
%\author[1]{Gy. Lipták} 
%\ead{lipgyorgy@scl.sztaki.hu}
%
%
%\author[4]{M. Kov\'acs}
%\ead{kovacs.mihaly@itk.ppke.hu}
%
%\author[4,1]{G. Szederk\'enyi}
%\ead{szederkenyi@itk.ppke.hu}
%
%
%
%\address[2]{Center for Geosciences and Geoengineering, Mines Paris - PSL University,
%	Fontainebleau, France.} 
%
%
%\address[3]{Department of Electrical Engineering,
%	Chalmers University of Technology,
%	G\"oteborg, Sweden.} 
%
%\address[1]{Institute for Computer Science and Control (SZTAKI),
%	Budapest, Hungary.} 
% 
%\address[4]{Faculty of Information Technology and Bionics,
%	Pázmány Péter Catholic University,
%	Budapest, Hungary.} 

%\nonumnote{The  authors  acknowledge  the  helpful  discussions  with  M. Vághy. Thank you for being our disciplined  colleague,  Pinar, some memories will never fade.    B.  Kulcsár  acknowledges the  contribution of Transport Area of Advance  at Chalmers  University  of Technology. The  project  has  been  partially  supported by European Comision(F-ENUAC-2022-0003) and Energimyndigheten (2023-00021) through the project E-LaaS:Energy optimal urban Logistics As A Service.  M. Kovács acknowledges the support of the grant NKFIH TKP2021-NVA-02. G. Szederkényi acknowledges the support of the projects NKFIH 145934 and RRF-2.3.1-21-2022-00006.}

\maketitle

% flatex input end: [title.tex]

%--------------------------------------------------------------------------------------
%--------------------------------------------------------------------------------------

%--------------------------------------------------------------------------------------------

\begin{abstract}	
	In this work, a family of finite volume discretization schemes for LWR-type first order traffic flow models (with possible on- and off-ramps) is proposed: the Traffic Reaction Model (TRM). These schemes yield systems of ODEs that are formally equivalent to the kinetic systems used to model chemical reaction networks. An in-depth numerical analysis of the TRM is performed. On the one hand, the analytical properties of the scheme (nonnegative, conservative, capacity-preserving, monotone) and its relation to  more traditional schemes for traffic flow models (Godunov, CTM) are presented. Finally, the link between the TRM and kinetic systems is exploited to offer a novel compartmental interpretation of traffic models. In particular, kinetic theory is used to derive  dynamical properties (namely persistence and Lyapunov stability) of the TRM for a specific road configuration. Two extensions of the proposed model, to networks and changing driving conditions, are also described.
\end{abstract}

%--------------------------------------------------------------------------------------------


\section{Introduction}

Macroscopic traffic flow models play an important role in the analysis, simulation, and control of traffic flows on networks \cite{SIRI2021,Ambrosio2009,Book_Treiber2013,Book_Papa,Papa2019,KARAFYLLIS2019228,coogan16,intro_traffic}. In this line of research,  first order traffic flow models, that describe the spatio-temporal evolution of the vehicular density through conservation laws,  predominate in the literature \cite{Piccoli2016,tordeux16}. These models, e.g. \cite{Piccoli2016,Book_Haight, lighthill_whitham,ctm_daganzo}, express the conservation of vehicles on the road by linking the density and the flux of vehicles in space and time.  In particular,  Lighthill and Whitham \cite{lighthill_whitham}, and Richards \cite{richards} suggested to further assume that the flux of vehicles can be expressed as function of the density called flux function (or fundamental diagram). The resulting traffic model, called Lighthill-Whitham-Richards (LWR) model, has been extensively studied and takes the form of a nonlinear hyperbolic Partial Differential Equation (hPDE) satisfied by the density of vehicles \cite{Bressan2009,Book_Treiber2013}. In general, no closed-form solution is available to solve such hPDEs and therefore \emph{appropriate} numerical approximations must be introduced.

Finite Volume Methods (FVMs) have been particularly praised for deriving numerical solutions of conservation laws, especially for their ability to accurately describe the singularities specific to these hPDEs (e.g. shock waves)  \cite{Book_Leveque1992,Bressan2009}. The Godunov scheme  \cite{godunov} is an instance of such methods widely used in traffic flow modeling for its advantageous theoretical and convergence properties \cite{Book_Leveque1992}. For instance, Daganzo proposes to equip the traditional Godunov scheme \cite{godunov} with static saturation functions (capacitated sending-receiving flows between cells), thus yielding the so-called Cell Transmission Model (CTM) \cite{ctm_daganzo}. This model is convenient in terms of describing capacity changes and has been the basis of many traffic network control solutions e.g. \cite{csikos2017}. 

%Other numerical schemes have also been implemented on LWR models, e.g. \cite{DAGANZO1995261,ZHANG2001337,LEO1992207}. 
The choice of numerical scheme is instrumental in deriving proper approximations of the solution of the considered hPDEs. First, the numerical approximations at hand need to preserve several important properties of the original hPDE like the nonnegativity of the solution or the fact that it should be bounded by the road capacity at all times. Besides, when it comes to FVMs, their stability \cite{FVM_book_2000} can be assessed by analyzing the systems of ordinary differential equations (ODEs) resulting from the space discretization of the PDE. These systems of ODEs may then have a connection to an already existing modeling framework thus facilitating their analysis. 

On the other hand, positive (or nonnegative) systems having the property that the state variables are always nonnegative have a distinguished role in systems and control theory due to the wide field of possible applications e.g., in thermodynamics, chemistry, biology or economy \cite{farina2000positive}. Moreover, the positivity of the states efficiently supports dynamical analysis and control design \cite{haddad2010nonnegative}. Within nonnegative systems, kinetic models (also called chemical reaction networks or CRNs) and the closely related compartmental models are good descriptors of complex nonlinear dynamics, and their mathematical and corresponding graph structure are intuitive and useful to state strong results on qualitative dynamics \cite{erdi1989mathematical,jacquez1993qualitative}. The theory of kinetic systems is a rapidly developing area of increasing interest where the most significant contributions are related to the existence and uniqueness of equilibria and robust stability even when several model parameters (possibly including time delays) are not precisely known or unknown \cite{Feinberg2019}. Therefore, it is of general interest to transform originally non-chemical models to kinetic form, and interpret or discover their properties from this point of view \cite{Samardzija1989}.

Building on these ideas, we propose and study a family of finite volume approximations of the LWR model, called Traffic Reaction Model (TRM), with direct links to the positive kinetic systems used to model chemical reaction networks \cite{Feinberg2019,SzeMaHan2018}. The benefits of this analogy are multiple. On the one hand, the resulting model provides a new interpretation of the discrete dynamics of the LWR model as a compartmental chemical system in which the concentration of \q{particles} of free and occupied space are exchanged between adjacent road cells. This analogy also allows to straightforwardly extend the LWR model to better represent traffic conditions that can occur in real life (e.g. merging, changing road conditions, traffic lights).

\changednew{Invariance and symmetry properties have extensively been researched in kinematic wave theory e.g. \cite{Newell1, Newell2,Newell3,Daganzo2005,LAVAL201317}. In the latter two papers, duality plays a key role to bridge the gap between coordinates in which traffic flow models are presented. A key consequence is that variational theory enables us to connect homogeneous first order traffic flow models with proper initial and boundary value problems. In \cite{LAVAL2016168}, the authors explore the symmetry of kinematic wave models under coordinates transformation and dual variable changes. When dualized, it is shown that travel time, delay, and flow remain invariant. In our paper, flow invariance is in the epicenter of TRM (through the re-definition of the flux by the dual state variable). The dual state variable is the free space concentration or dual density, which ensures the invariance of the nonlinear flux function %$f(\rho)=g(\rho,\nu)$ 
via a non-unique parametrization and decomposition. However, unlike in \cite{LAVAL2016168}, we do not try to explicitly transform the conservation laws and define new coordinates to ease delay or capacity calculations (for triangular fundamental diagrams). Instead, we use a nonlinear flux function, and decompose and parameterize it with the dual state variable (free space), where no coordinates transformation is required. Via the dual state variable and with the help of the flux decomposition, the numerical segmentation of the hyperbolic PDE (TRM) preserves several fundamental properties (e.g. non-negativity, capacity, etc). Furthermore, TRM can handle (under some mild conditions) inhomogeneity in PDEs and can accommodate a wide set of nonlinear flux functions. Finally, the system of nonlinear ODEs obtained by the TRM scheme can be interpreted as particular instances of kinetic dynamical systems which can be formally represented as chemical reactions. This special description opens the way to leverage the extensive literature on reaction networks to derive and analyze the dynamical properties of the system (such as robustness or stability) \cite{Feinberg1987,Chaves2002,Son:2001,Shinar2010}, and even propose novel optimal control strategies and traffic management techniques \cite{Shinar2010,7058380} towards \cite{CSIKOS2017120,DABIRI2015585,7539538}.} 

The main contributions of this work are as follows. Firstly, we propose an in-depth analysis of the TRM as a numerical scheme. We prove that it defines is consistent, monotone, nonnegative, capacity-preserving, and conservative numerical schemes. We also draw parallels between the TRM and existing methods, showing in particular that the TRM is in fact equivalent to a generic formulation of Daganzo's CTM, and that the Godunov scheme can be seen as particular instance of TRM. Consequently, both methods can also be treated as particular kinetic systems and the results from kinetic system theory can also be extended to them.

Secondly, we provide the arguments allowing to bridge the gap between numerical schemes for traffic models and kinetic systems for chemical reaction networks. We provide  detailed physical interpretations of how traffic flow is apprehended when seen as a kinetic system. We then leverage the link between the TRM and kinetic systems to provide an analysis of some dynamical properties of the TRM in a particular case. Persistence is introduced as new concept for traffic flow models based on this parallel. Also the structure and stability of the equilibria of the system are studied using Lyapunov theory for kinetic systems.
 
Finally, we present two direct extensions of the TRM that directly follow from the kinetic interpretation of traffic proposed in this work. The first one accounts for changing road conditions. Such extensions can be used to accurately model real traffic data set in state estimation and short-term prediction applications, as advocated in \cite{pereira2022parameter,pereira2022short}. The second extension allows us to describe traffic on road networks. \changednew{In this aspect, the TRM hints to model the link and node behaviour identically using  reaction rate governed mixing compartments, giving rise to a smooth interface between link and node modeling framework. Furthermore, the network extension of TRM fits in the genealogy of first order node models presented in \cite{TAMPERE2011289}. The network extended TRM fulfills the requirements for node models described in \cite{TAMPERE2011289}, where most of the properties are preserved thanks to the applied non-negative discretization scheme. Interestingly, many of the requirements are of continuous nature (e.g. conservation, capacity, non-negativity, or invariance) and no external constraints need to be enforced to guarantee them. It is important to note that the \emph{Network TRM} presented in this paper covers only non-signalized intersections.}

The outline of the paper is as follows. In \Cref{sec:trm}, we introduce the TRM as a finite volume scheme for the LWR model and describe its main properties. We then expose in \Cref{sec:kinetic} the coupling that exists between kinetic reaction theory and the TRM. In \Cref{sec:aTRN}, we provide an analysis of the TRM through its interpretation as particular case of CTM  and a numerical experiment testing the convergence and accuracy of the TRM. In \Cref{sec:dyn}, we derive some dynamical properties (namely persistence and Lyapunov stability) of the TRM when considered on a circular road. Finally, we present in \Cref{sec:ext} the two extensions of the TRM dicussed above.

% flatex input end: [Content/intro.tex]

%--------------------------------------------------------------------------------------%--------------------------------------------------------------------------------------

\section{Nonnegative discretization of hyperbolic PDEs}\label{sec:trm}

% flatex input: [Content/nonnegative.tex]

The Lighthill--Whitham--Richards (LWR) traffic model \cite{lighthill_whitham} is the first-order macroscopic traffic flow model defined by the following partial differential equation (PDE) on the domain $\mathbb{R}\times \overline{\mathbb{R}}_+$
\begin{equation}
\label{eq:scalar_consv_pde}
\partial_t\rho+\partial_xf(\rho)_x = r - s,
\end{equation}
This PDE models the space-time evolution of a conserved quantity $\rho(x, t) \in \Omega=\left[0, \rho_{\max}\right]$ corresponding to to the traffic density on a road, where $\rho_{\max}$ denotes the maximal vehicle density.  The flux of vehicles is represented $f\circ\rho$, for some  continuous and at least once differentiable $f : \mathbb{R} \to \mathbb{R}$  called flux function. In particular, we denote by $\rho_{\text{crit}}$  the critical density value at which the flux attains its maximal value $f_{\max}$  (i.e. $f_{\text{max}}=f(\rho_{\text{crit}})$). Finally, the functions $r, s$ model the on- and off-ramps, respectively.  

\begin{remark}
    We consider physically-relevant (weak) solutions $\rho$ of PDE~\eqref{eq:scalar_consv_pde} also called %mild solutions or %
entropy solutions. We refer the reader to Appendix \ref{app:3} for a result justifying the existence and uniqueness of an entropy solution, and the regularity of this solution under non-zero sink and source terms.
\end{remark}

To mirror the vehicle density $\rho$, we introduce a function $\nu$ modeling the density of free space (i.e. unoccupied by a vehicle) along the road, and given by
\begin{equation}
    \nu=\rho_{\max} - \rho.
\end{equation}
We set $\nu_{\text{max}}=\rho_{\max}$. Note then that the free space density $\nu$ can be seen as a dual or adjoint variable to the density $\rho$. \footnote{\changednew{$\nu$ is the concentration or density of free spaces}}
In this work, we consider the case where the flux function $f$ has the form
\begin{equation}
\label{eq:decomposed_flux}
f(\rho) = g(\rho, \nu)=g(\rho, \rho_{\max}-\rho),
\end{equation}
where $g: \Omega\times \Omega \mapsto \overline{\mathbb{R}}_+$ satisfies the following assumptions\footnote{These requirements are inspired by the properties of chemical reaction rates, as it will be explained in \Cref{sec:kinetic}.}:
\begin{enumerate}[label=({A}{\arabic*})]
	\item\label{assum:1} $g$ is Lipschitz continuous w.r.t. both arguments, with associated Lipschitz constants $K_1 >0$ and $K_2>0$,
	\item\label{assum:2} $g$ is non-decreasing in each arguments,
	\item\label{assum:3} $g(\rho, 0) = g(0, \nu) = 0$ for all $\rho, \nu \in \Omega$ (which ensures that no vehicles is removed (resp. added) if the road is empty (resp. at capacity).
\end{enumerate}

\begin{Ex}
	Any flux function $f$ that can be written as
	\begin{equation*}
	f(\rho)=g_1(\rho)g_2(\rho_{\max}-\rho), \quad \rho\in \Omega,
	\end{equation*}
	where $g_1$ and $g_2$ are non-decreasing Lipschitz-continuous functions  such that $g_1(0)=g_2(0)=0$ satisfy Assumptions \ref{assum:1}-\ref{assum:3}.  In particular, if we take $g_1(\rho)=\rho$, then the function $\rho\mapsto g_2(\rho_{\max}-\rho)$ can be interpreted as the speed-density relationship of the fundamental diagram. Within this framework, taking $g_2(\nu)\propto \nu$ yields the Greenshields (quadratic) flux function (cf. \Cref{fig:trap_fd}).
	We can also retrieve trapezoidal fundamental diagrams by considering $g_2$ as 
	\begin{equation*}
	g_2(\nu)=
	\begin{cases}
	v_{\max} \frac{\rho_{\max}-\nu_1}{\rho_{\max}-\nu}\frac{\nu}{\nu_2} & \text{if } 0 \le \nu \le \nu_2,\\
	v_{\max} \frac{\rho_{\max}-\nu_1}{\rho_{\max}-\nu} & \text{if } \nu_2 \le \nu \le \nu_1, \\
	v_{\max} & \text{if } \nu_1 \le \nu \le \nu_{\max}, \\
	\end{cases}
	\end{equation*}
	where $v_{\max}$ denotes the free flow speed, and $\rho_1=\rho_{\max}-\nu_1$, $\rho_2=\rho_{\max}-\nu_2$ denote the critical densities of the fundamental diagram which satisfy $0<\rho_1 \le \rho_2 < \rho_{\max}$.  	
	\Cref{fig:trap_fd} gives an example of a trapezoidal flux function obtained using this decomposition.
	
	\begin{figure}[t]
		\centering
		\begin{subfigure}[t]{0.4\textwidth}
			\includegraphics[width=\textwidth]{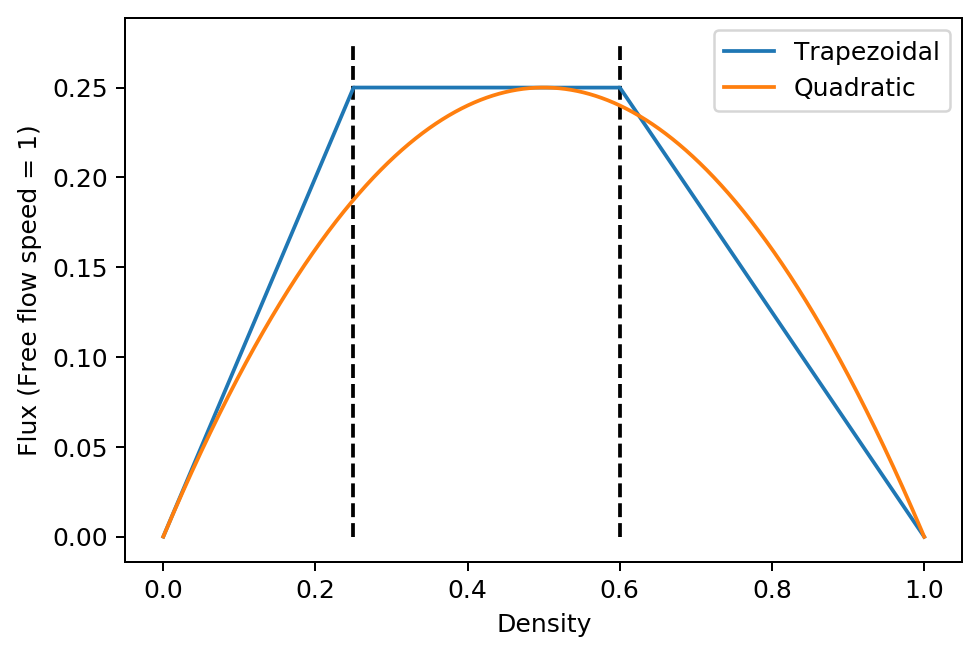}
			\vspace{-2em}
			\caption{Flux-density relationship}
		\end{subfigure}
	\hspace{2em}
		\begin{subfigure}[t]{0.4\textwidth}
			\includegraphics[width=\textwidth]{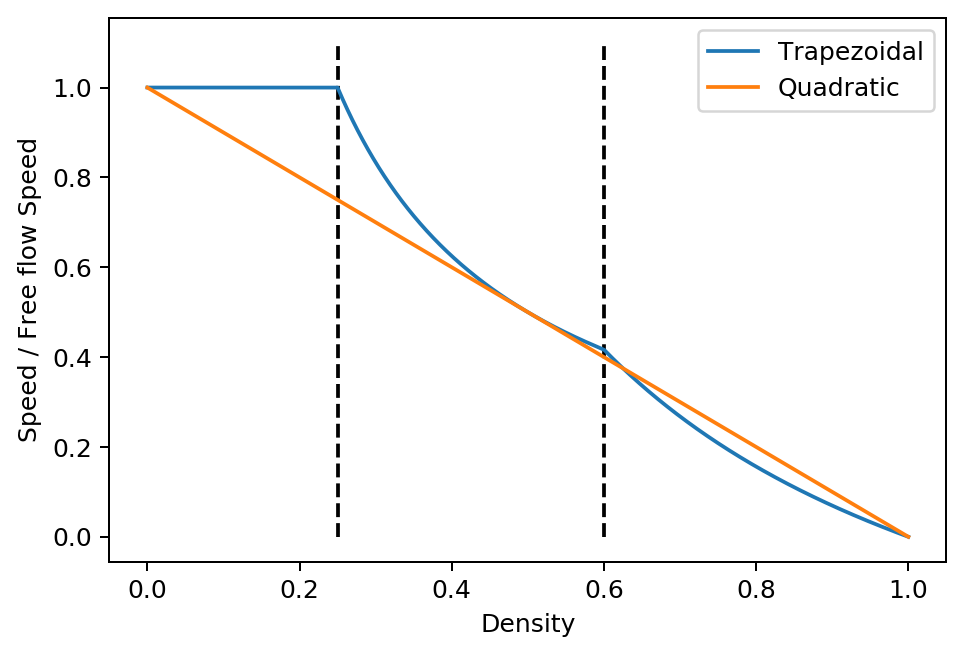}
			\vspace{-2em}
			\caption{Speed-density relationship (Function $\rho\mapsto g_2(\rho_{\max}-\rho)$)}
		\end{subfigure}
		\caption{Comparison between quadratic and trapezoidal fundamental diagrams. The maximal density $\rho_{\max}$ and the free flow speed $v_{\max}$ are taken equal to $1$. The black dotted lines are placed at the critical density values $\rho_1=0.25$ and $\rho_2=0.6$.}
		\label{fig:trap_fd}
	\end{figure}
\end{Ex}

Regarding the source and sink terms, we assume that they can be written for any $(x,\rho,t)\in\R\times\Omega\times\R_+$ as:
\begin{equation}
\label{eq:decomposed_source_sink}
    r(x, t, \rho) = 1_{\text{on}}(x) g_{\text{on}}(\rho_{\text{on}}(x,t), \nu),  \quad
    s(x, t, \rho) = 1_{\text{off}}(x) g_{\text{off}}(\rho, \nu_{\text{off}}(x,t)),
\end{equation}
where 
\begin{itemize}
	\item the functions $\rho_{\text{on}}$ and $\nu_{\text{off}}$ are the traffic density of the on-ramp and the free space density of the off-ramp, respectively, and are assumed to be 
	%piecewise Lipschitz functions (with a finite number of discontinuities) 
	taking values in $\Omega$;
	\item the functions $g_{\text{on}}:\Omega\times \Omega \mapsto \overline{\mathbb{R}}_+$  and $g_{\text{off}}:\Omega \times \Omega \mapsto \overline{\mathbb{R}}_+$ are the traffic flows of the on- and off-ramp, respectively, and are assumed to satisfy Assumptions \ref{assum:1}-\ref{assum:3};
	\item  the spatial position of the on- and off- ramp is described by indicator functions $1_{\text{on}}$ and $1_{\text{off}}$ defined as
	\begin{equation*}
	1_{\text{on}}(x)=\left\{ \begin{array}{lll}
	1, & \text{if } \underline{x}_{\text{on}} \le x \le \overline{x}_{\text{on}} \\
	0, & \text{otherwise} &
	\end{array}\right.,
    \quad
	1_{\text{off}}(x)=\left\{ \begin{array}{lll}
	1, & \text{if } \underline{x}_{\text{off}} \le x \le \overline{x}_{\text{off}} \\
	0, & \text{otherwise} &
	\end{array}\right.
	,
	\end{equation*}
	for some $\underline{x}_{\text{on}} \le \overline{x}_{\text{on}}$ and $\underline{x}_{\text{off}} \le \overline{x}_{\text{off}} $.
\end{itemize}

We use the Finite Volume Method (FVM) to spatially discretize (i.e. semi-discretize) the conservation law model \eqref{eq:scalar_consv_pde} with flux \eqref{eq:decomposed_flux}, source and sink \eqref{eq:decomposed_source_sink} functions. To do so, we start by slicing the road into cells of length $\Delta x>0$, and centered at the points $x_{i}= i\Delta x$, for $i\in\mathbb{Z}$. We denote by $I$ the set of road cell indices (hence $I=\mathbb{Z}$ for an infinite road) and denote by $x_{i-1/2}$ and $x_{i+1/2}$ the left and right boundary of the $i$-th cell, respectively, and by $\mathcal{C}_i=[x_{i-1/2}, x_{i+1/2}]$ the $i$-th cell. Then, the FVM scheme can be written as 
\begin{equation}
\label{eq:gen_scheme}
\dot{\rho}_i(t) = \frac{1}{\Delta x} \big(F(\rho_{i-1}, \rho_{i}) - F(\rho_{i}, \rho_{i+1}) 
+ R_{i}(\rho_i, t) - S_{i}(\rho_i, t)\big), \quad i \in I, \quad t\ge 0,
\end{equation}
where $\rho_i$ is an approximation of the average traffic density over the $i$-th cell,  the so-called numerical flux $F$ is given by
$F(u, v) = g(u, \rho_{\max} -v)$,
and the numerical source $R_i$ and sink $S_i$ terms are given (for $t\ge 0, i\in I$) by
\begin{equation*}
R_{i}(\rho, t) =\int_{\mathcal{C}_i} r(x,\rho,t)dx, \quad
S_{i}(\rho, t) = \int_{\mathcal{C}_i} s(x,\rho,t)dx.
\end{equation*}
We call \textit{Traffic Reaction Model} (TRM) the resulting model which describes traffic along a (discretized) road. Note that this model is control-oriented since changing $g$ and/or $g_{\text{on}}, g_{\text{off}}$ corresponds to speed limit control and/or ramp metering.
%where $1_{\text{on},i}=\int_{x_{i-1/2}}^{x_{i+1/2}}1_{\text{on}}(x)dx$, $1_{\text{off},i}= \int_{x_{i-1/2}}^{x_{i+1/2}}1_{\text{off}}(x)dx$. 
We also take as initial condition:
\begin{equation*}
\rho_i(0)=\frac{1}{\Delta x}\int_{\mathcal{C}_i}\rho(x,0)\,dx \in \Omega, \quad i \in I.
\end{equation*}
As such, the TRM defines a particular instance of finite volume scheme for scalar conservation laws \cite{leveque2002finite,chainais2001finite}.

Finally, note that in practice, finite roads should be considered. To be able to use the TRM to model traffic on a finite road, boundary conditions describing the incoming and outgoing traffic on the road have to be defined. Assuming that the road is now discretized into $P$ cells with indices $i\in I=\lbrace 1, \dots, P\rbrace$, the density of vehicles on these cells is modeled using the system of ODEs defined by~\eqref{eq:gen_scheme} while assuming that the quantities $\rho_0$ and $\rho_{P+1}$ are some specific time-dependent functions taking values in $\Omega$. For instance, these functions can be set so that the quantities $F(\rho_0, \rho_1)$ and $F(\rho_P, \rho_{P+1})$ match some known incoming and outgoing fluxes in the road. Another possibility is to set $\rho_{0}=\rho_P$ and $\rho_{P+1}=\rho_1$, which corresponds to periodic boundary conditions, or equivalently a circular road.

We now present a few numerical properties of the TRM, which hold whether the road is finite or not (cf. Appendix~\ref{sec:proofThm} for proof).

\begin{theorem}\label{thm:genFVM}
	The numerical flux $F$ 
	is consistent with the flux function $f$ (i.e. it satisfies for any $u\in\Omega$, $F(u,u)=f(u)$)  and is monotone. As for the TRM defined in~\eqref{eq:gen_scheme}, it preserves nonnegativity and capacity, and it is conservative finite volume scheme.
\end{theorem}

	It is important to note that, for a given $f$ the choice of $g$ is not unique. Hence, the TRM defines a family of numerical discretization schemes by means of the decomposition of the flux $f(\rho)$ in \eqref{eq:decomposed_flux}, which all share the properties described in Theorem \ref{thm:genFVM}.
	
\begin{Ex}\label{ex:g}\label{rem:god}
	Consider the so-called Greenshield flux function defined by
\begin{equation*}
    f(\rho) = \rho\, v(\rho), ~\text{with}~ v(\rho) = \omega(\rho_{\max} - \rho), \quad \rho\in \Omega,
\end{equation*}
where $v$ is the average vehicle speed, and $\omega = v_{\max} / \rho_{\max}$ with the free flow speed $v_{\max}>0$. 
Different decompositions of $f$ as in \eqref{eq:decomposed_flux} can be proposed, among which:
\begin{align}\label{eq:g1}
&g_\Mk(\rho, \nu) = \omega \rho \nu, \\
\label{eq:g2}
&g_\Gd(\rho, \nu) = \min(D(\rho), Q(\rho_{\max} - \nu)), \\
\label{eq:g3}
&g_\Cp(\rho, \nu) = D(\rho) Q(\rho_{\max} - \nu)/f_{\max},
\end{align}
%\begin{equation}\label{eq:g1}
%    g_\Mk(\rho, \nu) = \omega \rho \nu,
%\end{equation}
%\begin{equation}\label{eq:g2}
%    g_\Gd(\rho, \nu) = \min(D(\rho), Q(\rho_{\max} - \nu)),
%\end{equation}
%and
%\begin{equation}\label{eq:g3}
%    g_\Cp(\rho, \nu) = D(\rho) Q(\rho_{\max} - \nu)/f_{\max},
%\end{equation}
where $D(\rho)=f(\min\lbrace\rho, \rho_{\text{crit}}\rbrace)$, $Q(\rho) =f(\max(\rho, \rho_{\text{crit}}))$ are the supply and demand functions, respectively, and $f_{\max}$ denotes the maximal value of the flux. Note in particular that  the decomposition $g_\Gd(\rho,\nu)$ in \eqref{eq:g2} is the Godunov flux, meaning that the Godunov scheme can be seen as a particular instance of TRM.

\noindent In the remainder of this text, we refer to the decomposition $g_\Mk$ in \eqref{eq:g1} as a mass action kinetic (MAK) decomposition.
\end{Ex}

We further discretize the TRM by considering a time step $\Delta t>0$ and taking $t_k=k\Delta t$, $k\in\mathbb{N}_0$. We then define the fully-discrete TRM to be sequence of values $\lbrace \rho_i^k : i\in I, k\in\mathbb{N}_0\rbrace$ defined by the recurrence relations:
\begin{align}\label{eq:discrete_kinetic}
&\rho_i^{k+1}=\rho_i^{k}+\frac{\Delta t}{\Delta x}\big(F(\rho_{i-1}^{k},\rho_{i}^{k})-F(\rho_{i}^{k},\rho_{i+1}^{k})+q_i^k\big),\\
&\rho_i^0 =\rho_i(0) \in \Omega,
\end{align}
where the term $q_i^k$ is given by
\begin{equation*}
\begin{aligned}
q_i^k
=\frac{1}{\Delta t}\bigg[
&\int_{t_k}^{t_{k+1}}R_{i}(\rho_i^k, t) dt  -\int_{t_k}^{t_{k+1}} S_{i}(\rho_i^k, t) dt\bigg] .
\end{aligned}
\end{equation*}
In particular, if the road is finite, the quantities $\rho_j^k$, $j\in\lbrace 0, P+1\rbrace$ are defined from the functions $\rho_j$ as $\rho_j^k=(1/\Delta t)\int_{t_k}^{t_{k+1}}\rho_j(t)dt$.
The next result exposes the conditions under which fully-discrete TRM converges to the solution of the PDE. It relies on the following assumption on the source $r$ and sink $s$ terms.

\begin{assum}\label{assum:source}
	The source and sink term~$q=r-s$ of PDE \eqref{eq:scalar_consv_pde} is bounded and such that for any $\rho\in\R$,
	$(x,t)\mapsto q(x,t,\rho)$ is
	Lipschitz-continuous (with a constant independent of $u$) and for any $(x,t)\in\R\times\R_+$,
	$\rho\mapsto q(x,t,\rho)$ is
	locally Lipschitz-continuous (with a constant independent of $(x,t)$ and bounded $K_1$ and $K_2$ as defined in \ref{assum:1}).
\end{assum}

\begin{theorem} \label{thm:genFVM2}
	Let us assume that the  source and sink term~$q=r-s$ of PDE \eqref{eq:scalar_consv_pde} satisfy Assumption~\ref{assum:source} and that the following Courant-Friedrichs-Levy (CFL) condition is fulfilled:
	\begin{equation}
	\label{eq:CFL}
	\frac{\Delta t}{\Delta x}\leq \frac{1}{K_1+K_2} ,
	\end{equation}
	where $K_1$ and $K_2$ are defined in~\ref{assum:1}.
	%\cite{FVM_book_2000}.
 Then, the fully discrete TRM converges (in $L^1_{\text{loc}}(\R\times\R_+)$) towards the entropy solution of PDE \eqref{eq:scalar_consv_pde} as $\Delta x, \Delta t \rightarrow 0$, with $\Delta t / \Delta x$ kept fixed and satisfying the CFL condition. 
\end{theorem}

\begin{proof}
	This result is a direct application of Theorem 1 of \cite{chainais2001finite}, which holds since the discretized scheme is monotone (cf. \Cref{thm:genFVM}) and following from the assumptions on the regularity of the source and sink term $q$.
\end{proof}

% \begin{Ex}\label{rem:CFL2}
% 	Let us consider the flux function described in \Cref{ex:g}, and the TRM obtained using any choice of the decompositions \eqref{eq:g1}-\eqref{eq:g3}. Assuming that the source and sink terms satisfy the conditions in \Cref{thm:genFVM2}, then a sufficient condition for the (fully-discrete) TRM to converge to the entropy solution of PDE \eqref{eq:scalar_consv_pde} is that the time and space steps satisfy the CFL condition given by
% 	\begin{equation}
% 	\label{eq:CFL2}
% 	\frac{\Delta t}{\Delta x}\leq \frac{1}{2v_{\text{max}}}.
% 	\end{equation}
% \end{Ex}

\begin{remark}
In the absence of source and sink terms (and with an infinite road), but under the CFL condition \eqref{eq:CFL}, the discretized numerical scheme in \eqref{eq:discrete_kinetic} is %$\mathrm{i})$
$L^{\infty}$-stable, meaning that for any $(i, k)\in\mathbb{Z}\times \mathbb{N}_0$, $\rho_i^k \in [\min_{j\in\mathbb{Z}} \rho_j^0, \max_{j\in\mathbb{Z}}  \rho_j^0]$  (cf. Lemma $21.1$ in \cite{FVM_book_2000}). Consequently, the discrete numerical scheme preserves nonnegativity and capacity (as long as the initial condition is in $\Omega$). Furthermore, since the scheme is monotone and $L^{\infty}$-stable, it is Total Variation Diminishing (TVD) \cite[Theorems 15.4 \& 15.5]{Book_Leveque1992}. As a consequence, the TRM is in particular shock-capturing. It implicitly incorporates correct jump conditions (without oscillations) near discontinuities when $\Delta x$ and $\Delta t$ tend to $0$. Besides, the convergence of the (fully-discrete) TRM towards the PDE solution is of order at most 1 \cite{leveque2002finite}. We retrieved this behavior in numerical experiments exposed in Appendix~\ref{sec:numTRM}.
\end{remark}

% flatex input end: [Content/nonnegative.tex]

%--------------------------------------------------------------------------------------
%--------------------------------------------------------------------------------------

\section{Kinetic and compartmental description} \label{sec:kinetic}

In this section, we will show that the discretized traffic flow model \eqref{eq:gen_scheme} (on a finite road) is formally kinetic with appropriate assumptions on $g$, and it can be interpreted in a compartmental context.

\subsection{Kinetic models}\label{subsec:kinetic}
Here we characterize kinetic system models based on the notations used in \cite{Feinberg2019}, where more details can be found on this model class. We assume that a kinetic model contains $M$ species denoted by $\mathcal{X}_1,\dots,\mathcal{X}_M$, and the species vector is $\mathcal{X}=[\mathcal{X}_1~\dots~\mathcal{X}_M]^T$. The basic building blocks of kinetic systems are \textit{elementary reaction steps} of the form
 \begin{align}
 C_j \rightarrow C'_j,~j=1,\dots,R \label{eq:Reactions}
 \end{align}
where for $j=1,\dots,R$, $C_j=y_j^T \mathcal{X}$ and $C'_j={y'_j}^T \mathcal{X}$ are the \textit{complexes} and  $y_j,y'_j\in {\mathbb{N}}_{0}^M$ are integer coefficients. The transformation shown in~\eqref{eq:Reactions} means that during an elementary reaction step between the \textit{reactant complex} $C_j$ and \textit{product complex} $C_j'$, $[y_j]_i$ items (molecules) of species $\mathcal{X}_i$ are consumed, and $[y'_j]_i$ items of $\mathcal{X}_i$ are produced for $i=1,\dots,M$. The reaction \eqref{eq:Reactions} is an \textit{input} (resp. \textit{output}) \textit{reaction of species $\mathcal{X}_i$} if $[y_j']_i>0$ (resp. $[y_j]_i>0$). 

Let $\chi(t)\in\overline{\mathbb{R}}^M_{+}$ denote the state vector corresponding to $\mathcal{X}$ for any $t\ge 0$ (in a chemical context, $\chi$ is the vector of concentrations of the species in $\mathcal{X}$). Then the ODEs describing the evolution of $\chi$ in the kinetic system containing the reactions \eqref{eq:Reactions} are given by
\begin{align}\label{eq:general_kinetic}
\dot{\chi}=\sum_{i=1}^R \RR_i(\chi,t)[y_i' - y_i],~\chi(0)\in\overline{\mathbb{R}}^M_+    
\end{align}
where $\RR_i:\overline{\mathbb{R}}_+^M \times \mathbb{R}  \longrightarrow \mathbb{R}_{+}$ is the \textit{rate function} corresponding to reaction step $i$, and determines the velocity of the transformation. For the rate functions, we assume the following for $i=1,\dots,R$:
\begin{enumerate}[label=({A}{\arabic*})]
\setcounter{enumi}{3}
\item $\mathcal{K}_i(\cdot,t)$ is locally Lipschitz continuous, and $\mathcal{K}_i(\chi, \cdot)$ is piecewise locally Lipschitz continuous with a finite number of discontinuities, 
\item for $j=1,\dots,m$, $\mathcal{K}_i(\cdot, t)$ depends on $\chi_j$ if and only if $y_j>0$,
\item $\mathcal{K}_i\ge 0$, and for $j=1,\dots,m$, $\mathcal{K}_i(\chi,t) = 0$ whenever $\chi_j=0$ and $y_j>0$,
\item There exist continuous nonnegative and strictly monotone (w.r.t. a partial order on $\overline{\mathbb{R}}_+^M$) functions $\underline{\mathcal{K}}_i$, $\overline{\mathcal{K}}_i: \overline{\mathbb{R}}_+^M \longrightarrow \mathbb{R}_{\ge 0}$ such that $\underline{\mathcal{K}}_i(\chi) \le \mathcal{K}_i(\chi, t) \le \overline{\mathcal{K}}_i(\chi) $ for all $t\ge 0$
\end{enumerate}
These properties ensure the local existence and uniqueness of the solutions as well as the invariance of the nonnegative orthant for the dynamics in~\eqref{eq:general_kinetic}. 
From now on, a reaction from complex $C_i$ to complex $C_i'$ with rate function $\RR_i$ will be denoted as $\displaystyle C_i \rrarrow{\RR_i} C_i'$. We will also suppress the $t$ argument in the rate function if it does not explicitly depend on time.

A set of nonlinear ODEs given as $\dot{\chi}=f(\chi)$
is called \textit{kinetic} if it can be written in the form \eqref{eq:general_kinetic} with appropriate rate functions $\mathcal{K}_i$. We remark that the representation $\eqref{eq:general_kinetic}$ of a kinetic ODE is generally non-unique even if the rate functions are a priori fixed \cite{acs2016computing}.  

\subsection{Kinetic representation of the traffic flow model}

Using the above notions, we can give a kinetic interpretation for the TRM, where the species and the reaction steps have a clear physical meaning. For this purpose, we consider the TRM~\eqref{eq:gen_scheme} on a finite finite road (with $P$ cells).

\begin{figure}
	\centering
	\resizebox{0.75\textwidth}{!}{
		% flatex input: [reaction_rates.tex]

\begin{tikzpicture}[x=\textwidth,y=0.25\textwidth,scale=1, every node/.style={scale=1.25}]

\draw[-{Latex[length=3mm,width=2mm]},black,line width=0.35mm] (0,1) -- (1,1) node[anchor=west,black]{$x$};
\draw[black,line width=0.35mm] (0.2,1.05) node[anchor=south,black]{$x_{j-1}$} -- (0.2,0.95) ;
\draw[black,line width=0.35mm] (0.5,1.05) node[anchor=south,black] (tot) {$x_{j}$} -- (0.5,0.95) ;
\draw[black,line width=0.35mm] (0.8,1.05) node[anchor=south,black]{$x_{j+1}$} -- (0.8,0.95) ;

\draw[line width=0.5mm,black] (0.01,0.75) -- (0.99,0.75);
\draw[line width=0.5mm,black](0.01,0.45) -- (0.99,0.45);

\draw[line width=0.5mm,black] (0.01,0.15) -- (0.99,0.15);
\draw[line width=0.5mm,black](0.01,-0.5) -- (0.99,-0.5);

\draw[dashed,red] (0.05,1)  -- (0.05,-0.6);
\draw[dashed,red] (0.35,1)    -- (0.35,-0.6) node[anchor=north] {{ $\displaystyle S_{j}+N_{j-1}\mathop{\rightarrow} N_j + S_{j-1}$}};
\draw[dashed,red] (0.65,1)  -- (0.65,-0.6) node[anchor=north] {{ $S_{j+1}+N_{j}\mathop{\rightarrow} N_{j+1} + S_{j}$}};
\draw[dashed,red] (0.95,1)  -- (0.95,-0.6);

\draw[red,line width=0.5mm] (0.05,0.15)  -- (0.05,-0.5);
\draw[red,line width=0.5mm] (0.35,0.15)  -- (0.35,-0.5);
\draw[red,line width=0.5mm] (0.65,0.15)  -- (0.65,-0.5);
\draw[red,line width=0.5mm] (0.95,0.15)  -- (0.95,-0.5);

\draw[red,line width=0.5mm] (0.05,0.75)  -- (0.05,0.45);
\draw[red,line width=0.5mm] (0.35,0.75)  -- (0.35,0.45);
\draw[red,line width=0.5mm] (0.65,0.75)  -- (0.65,0.45);
\draw[red,line width=0.5mm] (0.95,0.75)  -- (0.95,0.45);

\draw (0.2,0.25) node [black] {$j-1$};
\draw (0.5,0.25) node [black] {$j$};
\draw (0.8,0.25) node [black] {$j+1$};

\draw (0.825,0.6) node {\includegraphics[scale=0.7]{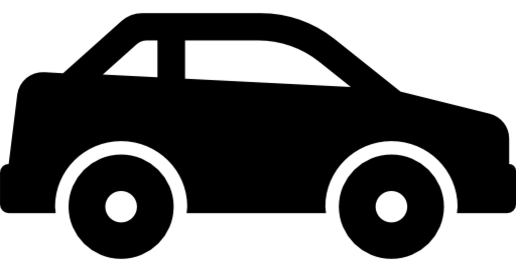}};
\draw (0.1,0.6) node {\includegraphics[scale=0.7]{car.png}};
\draw (0.2,0.6) node {\includegraphics[scale=0.7]{car.png}};
\draw (0.4,0.6) node {\includegraphics[scale=0.7]{car.png}};
\draw (0.5,0.6) node {\includegraphics[scale=0.7]{car.png}};
\draw (0.575,0.6) node {\includegraphics[scale=0.7]{car.png}};

\draw (0.5,-0.25) node [circle,draw=gray,inner sep=0cm,minimum size=7.5mm] {$N_j$};
\draw (0.46,0.01) node [circle,draw=gray,inner sep=0cm,minimum size=7.5mm] {$N_j$};
\draw (0.57,-0.3) node [circle,draw=gray,inner sep=0cm,minimum size=7.5mm] {$S_j$};
\draw (0.39,-0.01) node [circle,draw=gray,inner sep=0cm,minimum size=7.5mm] {$S_j$};
\draw (0.61,-0.03) node [circle,draw=gray,inner sep=0cm,minimum size=7.5mm] {$N_j$};
\draw (0.39,-0.3) node [circle,draw=gray,inner sep=0cm,minimum size=7.5mm] {$N_j$};
\draw (0.54,-0.02) node [circle,draw=gray,inner sep=0cm,minimum size=7.5mm] {$N_j$};

\draw (0.085,0) node [circle,draw=gray,inner sep=0cm,minimum size=7.5mm] {$S_{j-1}$};
\draw (0.15,0.0) node [circle,draw=gray,inner sep=0cm,minimum size=7.5mm] {$S_{j-1}$};
\draw (0.2,-0.2) node [circle,draw=gray,inner sep=0cm,minimum size=7.5mm] {$N_{j-1}$};
\draw (0.3,-0.31) node [circle,draw=gray,inner sep=0cm,minimum size=7.5mm] {$S_{j-1}$};
\draw (0.25,-0.03) node [circle,draw=gray,inner sep=0cm,minimum size=7.5mm] {$S_{j-1}$};
\draw (0.115,-0.27) node [circle,draw=gray,inner sep=0cm,minimum size=7.5mm] {$N_{j-1}$};
\draw (0.315,0.0) node [circle,draw=gray,inner sep=0cm,minimum size=7.5mm] {$N_{j-1}$};

\draw (0.9,-0.3) node [circle,draw=gray,inner sep=0cm,minimum size=7.5mm] {$S_{j+1}$};
\draw (0.785,0.0) node [circle,draw=gray,inner sep=0cm,minimum size=7.5mm] {$S_{j+1}$};
\draw (0.85,-0.02) node [circle,draw=gray,inner sep=0cm,minimum size=7.5mm] {$S_{j+1}$};
\draw (0.72,0.015) node [circle,draw=gray,inner sep=0cm,minimum size=7.5mm] {$S_{j+1}$};
\draw (0.71,-0.25) node [circle,draw=gray,inner sep=0cm,minimum size=7.5mm] {$N_{j+1}$};
\draw (0.82,-0.3) node [circle,draw=gray,inner sep=0cm,minimum size=7.5mm] {$S_{j+1}$};
\draw (0.915,0.015) node [circle,draw=gray,inner sep=0cm,minimum size=7.5mm] {$N_{j+1}$};

%\draw[latex-latex] (0.35,0.0) -- (0.65,0.0) node[midway, below] {$\Delta x$};

\end{tikzpicture}

% flatex input end: [reaction_rates.tex]

	}
	\caption{Representation of the Traffic Reaction Model interpretation of traffic flow. The discretized road (Top row) is seen as a sequence of compartments (Bottom row) containing \q{molecules} of free space $S$ and occupied space $N$ (represented as circles). The flow of vehicles along the road is then modeled by chemical reactions between the compartments (written in red).}
	\label{fig:rr}
\end{figure}

Let us model each road cell $i\in I=\lbrace 1, \dots, P\rbrace$ as a compartment containing two homogeneously distributed species: units of occupied space $N_i$ and units of free space $S_i$. The flow of vehicles between the consecutive cells (compartments) $i-1$ and $i$ (for $i\in\lbrace 2, \dots, P\rbrace$) is then represented by chemical reactions converting occupied space in one cell into occupied space on the next one, as follows
\begin{equation}
     N_{i-1} + S_i \rrarrow{\RR_{i-1,i}} N_i + S_{i-1}, \quad i\in\lbrace 2, \dots, P\rbrace, \label{eq:VehTrans1}
\end{equation}
where $\RR_{1-1,i}$ is the corresponding reaction rate.
This model is represented in \Cref{fig:rr}.
Note that, in this representation, units of space can flow from cell $i-1$ to cell $i$ provided that a unit of occupied space is present in cell $i-1$ and that a unit of free space is present in cell $i$. By recalling condition (A6) in Subsection \ref{subsec:kinetic}, it is easy to see that the model in Eq. \eqref{eq:VehTrans1}naturally enforces that there can be no flow from cell $i-1$ to cell $i$ if cell $i-1$ is empty (i.e. there is no unit of occupied space $N_{i-1}$) or if cell $i$ is full (i.e. there is no unit of free space $S_i$). 

Regarding the on and off-ramps in a given cell $i \in I$, they are respectively described by the following reaction steps
\begin{equation}\label{eq:VehTrans_on}
    S_i \rrarrow{\RR_{\text{on},i}} N_i , \quad
    N_i  \rrarrow{\RR_{\text{off},i}} S_i 
\end{equation}
where $\RR_{\text{on},i}$ and $\RR_{\text{off},i}$ are the reaction rates. These transformations  express that during an elementary step of on-ramp (resp. off-ramp) reaction, a unit of free space (resp. occupied space) is consumed, and turned into a unit of occupied space. As for the boundary compartments $j\in\lbrace 1, P\rbrace$, they are characterized by similar reaction steps, namely
\begin{equation}
    S_1 \rrarrow{\RR_{\text{in}}} N_1, \quad N_{P} \rrarrow{\RR_{\text{out}}} S_P 
    \label{eq:breac}
\end{equation}
where the reaction rates $\RR_{\text{in}}$ and $\RR_{\text{out}}$ are set to reflect the boundary conditions.

Let $n=(n_i : i\in I)$ and $s=(s_i : i \in I)$ denote the concentrations of occupied and free space in the compartment $i\in I$. Within our analogy, these quantities can respectively be seen as the equivalent of the density of vehicles and the density of free space in the cell. Then the kinetic differential system~\eqref{eq:general_kinetic} corresponding to the kinetic model with reaction steps \eqref{eq:VehTrans1}-\eqref{eq:VehTrans_on} is the system of ODEs defined by
\begin{align}
\label{eq:n_i_dot}
\dot{n}_i &=  \RR_{i-1,i}(n_{i-1},s_i) - \RR_{i,i+1}(n_i,s_{i+1})
 + \RR_{\text{on},i}(s_i,t) -\RR_{\text{off},i}(n_i,t), \quad i\in I,\\
 \label{eq:s_i_dot}
\dot{s}_i &=   \RR_{i,i+1}(n_i,s_{i+1})- \RR_{i-1,i}(n_{i-1},s_i)
+ \RR_{\text{off},i}(n_i,t) - \RR_{\text{on},i}(s_i,t) , \quad i\in I,
\end{align}
where we adopt the notation $\RR_{0,1}(n_0, n_1)\equiv \RR_{\text{in}}(n,t)$ and $\RR_{P,P+1}(n_P, n_{P+1})\equiv \RR_{\text{off}}(n,t)$.

 Note that $ \dot{n}_i+\dot{s}_i=0$. Hence, $c_i=n_i+s_i$ is a first integral (conserved quantity) of the dynamics which can be interpreted as the maximal vehicle density in cell $i$. 
 It is also clear that the conservation of $c_i$ guarantees the boundedness of the solutions of \eqref{eq:n_i_dot}-\eqref{eq:s_i_dot}.
Substituting then $s_i=c_i-n_i$ into~\eqref{eq:n_i_dot} (and assuming that $ 0 \leq n_i(0)\le c_i$) gives (for $i\in I$)
\begin{equation}
\dot{n}_i =~  \RR_{i-1,i}(n_{i-1}, c_i-n_i) - \RR_{i,i+1}(n_i, c_{i+1}-n_{i+1})
+ \RR_{\text{on},i}(c_i-n_i,t) -\RR_{\text{off},i}(n_i,t)
\label{eq:n_i_dot2}
\end{equation}

This last equation is equivalent to the TRM equation~\eqref{eq:gen_scheme} when taking $n_i=\rho_i$ to be the vehicle density of the $i$-th cell, $c_i=\rho_{\text{max}}$ to be the maximal cell capacity, and defining the reaction rates as
\begin{equation}\label{eq:rr}
\RR_{i,i+1}(\rho,\nu)=g(\rho,\nu)/\Delta x,\quad i\in \lbrace 0, \dots, P\rbrace,
\end{equation}
and for $i\in I$, $\RR_{\text{on},i}(\rho,t)=R_i(\rho,t)/\Delta x$, and $\RR_{\text{off},i}(\rho,t)=S_i(\rho,t)/\Delta x$. A notable case is when we consider TRM with the decomposition $g=g_\Mk$ introduced in~\Cref{ex:g}. Then, the resulting reaction rates are exactly those obtained when considering the kinetic system under the so-called mass action kinetics~\cite{Chellaboina2009} assumption.

More generally, the kinetic interpretation of the TRM provides a new insight into the definition of fundamental diagrams, as modeled by the function $\rho \mapsto f(\rho)$ linking the flux to the density of vehicles. Indeed, the decomposition $g$ used to define $f$, which involves both the density of vehicles $\rho$ and its dual the free space density $\nu$, has a clear physical and kinetic interpretation: it models how free space is turned into occupied space by the flow of vehicles. It is clearly visible from Eq. \eqref{eq:rr} that there is a simple proportional relationship between the kinetic reaction rates and the speed-density relationship of the fundamental diagrams. This gives a transparent and physically meaningful mapping between traffic flows and kinetic models.
Therefore, two main conclusions can be drawn. First, the TRM offers a complementary microscopic point-of-view on traffic where the particles being modeled are not vehicles but units of free and occupied space, and which has a clear and straightforward link to the quantities modeled by the macroscopic traffic model. Second, the TRM allows to propose and to interpret fundamental diagram relations based on a microscopic model for the transformations of free space and occupied space.

%\changednew{It must be remarked here that the potential importance of empty spaces (voids) in traffic dynamics has already been addressed in the literature. In \cite{laval2006lane} it was first shown that voids created by lane changing vehicles finally reduce the flow. The effect of voids was also taken into consideration to give a new model for capacity drops at merges in \cite{leclercq2011capacity}.
%However, in our model, empty spaces are only produced parallelly to the transition of vehicles from one cell to the other since we do not consider lane changes or spatial distribution of individual vehicles within a cell.}

%\textbf{Gabor some extensions!!!!}

\subsection{Application: Dynamical analysis of the ring topology} \label{sec:dyn}
% flatex input: [Content/stab.tex]

We conclude thus section with an example of how the kinetic interpretation of the TRM can be leveraged to deduce some properties of the model using kinetic system theory. In particular, we consider the special case when the TRM with $P$ compartments has a ring topology (i.e.  periodic boundary conditions) without on and off-ramps, and study the persistence of the dynamics and then the stability of the equilibria. 

\subsubsection{Persistence of the dynamics using kinetic theory}

Following the kinetic interpretation made in~\Cref{sec:kinetic}, the TRM can be analyzed using the theory of kinetic systems. Indeed, the TRM is equivalent to the kinetic system defined by the reactions~\eqref{eq:VehTrans1} and~\eqref{eq:breac}. Note in particular that the boundary reactions~\eqref{eq:breac} can be merged into a single reaction linking the first and last cells, and given by
\begin{equation*}
    N_{P} + S_1 \rrarrow{\RR_{P, 1}} N_1 + S_{P}.
\end{equation*}
Since the reaction rates $\mathcal{K}_{i,i+1}$ are given by~\eqref{eq:rr} and do not depend explicitly on time, and assuming that they are real analytic functions, we can use the results of \cite{angeli2007petri} on the persistence of kinetic systems.

A nonnegative kinetic system is called \textit{persistent} if no trajectory that starts in the interior of the orthant $\mathbb{R}^P_+$ has an $\omega$-limit
point on the boundary of $\mathbb{R}^P_+$.
A non-empty set of species $\Sigma\subset \mathcal{X}$ is called a \textit{siphon} if each input reaction associated to $\Sigma$ is
also an output reaction associated to $\Sigma$. A siphon is \textit{minimal} if it does not contain (strictly) any other siphons. Naturally, the union of siphons is a siphon.
The sufficient conditions for the persistence of the dynamics are the following \cite{angeli2007petri}: (1) there exists a positive linear conserved quantity (i.e., a first integral) for the dynamics, and (2) each siphon contains a subset of species, the state variables of which define a nonnegative linear first integral for the dynamics.

It is straightforward to see that Condition (1) is fulfilled, since $I_0=\sum_{i=1}^P s_i+n_i$ is a conserved quantity for the system. For Condition (2), it can be checked that the only minimal siphons in the ring topology are: $\Sigma_N=\{N_1,\dots,N_P \}$, $\Sigma_S=\{S_1,\dots,S_P\}$, and $\Sigma_i=\{N_i,S_i\}$ for $i=1,\dots,P$. Since $I_N=\sum_{i=1}^P n_i$, $I_S=\sum_{i=1}^P s_i$, and $I_i=n_i+s_i$ for $i=1,\dots,P$ are also linear first integrals containing the state variables of $\Sigma_N$, $\Sigma_S$, and $\Sigma_i$ respectively, Condition (2) is fulfilled. Therefore, the dynamics of the ring topology are persistent for any numerical flux $F$ for which the corresponding reaction rates $\mathcal{K}_{i,i+1}$ satisfy the mild conditions described above. Note that in order to prove the persistence of the dynamics for a wide class of flux functions just by exploiting the ring structure, we needed both sets of state variables representing the density of vehicles as well as the density of free space.

\subsubsection{Structure of equilibrium and stability}

We first analyze the possible equilibrium points of the system. The proof of the next result can be found in Appendix~\ref{appen:proof_eq}.

\begin{proposition}\label{prop:eq_pt}
For the system of ODEs~\eqref{eq:gen_scheme} resulting from the mass action kinetic TRM with a ring topology,  the only possible equilibrium point $(\rho_1^*, \dots, \rho_P^*)$ is the one satisfying for any  $k\in I$
\begin{equation*}
    \rho_k^* = \bar{\rho}=\frac{1}{P}\sum_{i=1}^P\rho_i(0)
\end{equation*}
\end{proposition}

Hence, as one could have expected, the only possible equilibrium point for the system is the one where the vehicles distribute themselves uniformly across the (circular) road. Clearly there exist infinitely many equilibria for the system, but there is one unique positive equilibrium within each equivalence class $\mathcal{H}_c=\{\rho\in\Omega~|~\sum_{i=1}^P\rho_i = c \}$ defined by the conservation of vehicles. The next result elaborates on the stability of this equilibrium.

\begin{proposition}\label{prop:entropyLyapunov}
The equilibrium $\bar \rho$ is stable on $\Omega$.
\end{proposition}

\begin{proof}
We use the entropy-like Lyapunov function candidate well-known from the theory of reaction networks \cite{Feinberg1987}:
\begin{equation}\label{eq:LogLyap}
V(\rho) = \sum_{i=1}^P \bigg(\rho_i \,\left[\log\left(\frac{\rho_i}{\overline{\rho}}\right) - 1\right] + \overline{\rho}\bigg).
\end{equation}
Then, 
\begin{equation*}
\begin{aligned}
\dot{V} &= \frac{\omega}{\Delta x}\sum_{i=1}^P \log\left(\frac{\rho_i}{\overline{\rho}}\right)  \left[\rho_{i-1} (\rho_{\max} - \rho_{i}) - \rho_i(\rho_{\max} - \rho_{i+1})\right],
\end{aligned}
\end{equation*}
which gives
\begin{equation*}
\begin{aligned}
\dot{V} 
=~ &\frac{\omega \overline{\rho}}{\Delta x}\sum_{i=1}^{P-1}\frac{\rho_i(\rho_{\max} - \rho_{i+1})}{\overline{\rho}} \left[\log\left(\frac{\rho_{i+1}}{\overline{\rho}}\right) - \log\left(\frac{\rho_i}{\overline{\rho}}\right)\right] 
 +\frac{\rho_P(\rho_{\max} - \rho_{1})}{\overline{\rho}} \left[\log\left(\frac{\rho_{1}}{\overline{\rho}}\right) - \log\left(\frac{\rho_P}{\overline{\rho}}\right)\right],
\end{aligned}
\end{equation*}
then we use the inequality $e^a(b-a)\leq e^b - e^a$ to give an upper bound such that
\begin{equation*}
\begin{aligned}
\frac{\Delta x}{\omega}\dot{V}  \leq 
&\sum_{i=1}^{P-1} (\rho_{\max} - \rho_{i+1})  \left(\rho_{i+1} - \rho_i\right) 
  + (\rho_{\max} - \rho_{1})   \left(\rho_{1} - \rho_P \right) \\
%rac{\Delta x}{\omega}\dot{V} \leq 
%\sum_{i=1}^{P-1}(\rho_{\max} - \rho_{i+1})  \left(\rho_{i+1} - \rho_i\right) + (\rho_{\max} - \rho_{1})   \left(\rho_{1} - \rho_P \right)\\
% &=\sum_{i=1}^{P-1}\rho_{\max}\,\rho_{i+1} - \rho_{\max}\,\rho_i - \rho_{i+1}^2 + \rho_{i+1}\rho_i  \\
% & + \rho_{\max}\,\rho_{1} - \rho_{\max}\,\rho_P - \rho_{1}^2 + \rho_{1}\rho_P  \\
& = \sum_{i=1}^{P-1} \rho_{i+1}\rho_i - \rho_{i+1}^2 + \rho_{1}\rho_P - \rho_{1}^2 
 = \sum_{i=1}^{P} -\frac{1}{2}(\rho_{i}^2 - 2\rho_{i+1}\rho_i + \rho_{i+1}^2) 
 =\sum_{i=1}^{P}  -\frac{1}{2}(\rho_{i} -\rho_{i+1})^2 \leq 0.
\end{aligned}
\end{equation*}
Note that the entropy-like Lyapunov function $V$ does not depend on the model parameters, and $\dot{V}<0$ outside the equilibria . We also remark that the theory of kinetic and compartmental systems can still be used when the spatial discretization is non-uniform along the ring. In such a case, the analytical computation of the equilibrium is not straightforward. However, we know that the dynamics is persistent, and for any $c>0$, there exists a unique strictly positive equlibrium in $\mathcal{H}_c$ which is stable with the Lyapunov function defined in Eq. \eqref{eq:LogLyap}, and asymptotically stable within $\mathcal{H}_c$.
\end{proof}
It is important to add that the special case studied here can be greatly generalized using recent results. Persistence and stability can be similarly proved not only for the MAK decomposition used here, but for very general reaction rates and thus for a really wide set of fundamental diagrams for arbitrary strongly connected networks \cite{szederkenyi2022persistence}. Moreover, we can handle explicitly time-varying transition rates (possibly resulting from changing conditions and/or external control) and determine a whole family of logarithmic Lyapunov functions in that case \cite{vaghy2023persistence}.

%Moreover, similar properties can be proved for more complex network topologies, and more general flux (transition rate) functions including the MAK decomposition used in this example. For the details, see \cite{szederkenyi2022persistence,vaghy2022lyapunov}.}

%\textbf{Gabor extensions of the stability results! Traffic control.}

% flatex input end: [Content/stab.tex]

%--------------------------------------------------------------------------------------
%--------------------------------------------------------------------------------------

\section{Analysis of TRM: Comparison with the CTM} \label{sec:aTRN}

In this section, we formally show that the TRM is equivalent to the Cell Transmission Model (CTM) \cite{ctm_daganzo} under some specific parametrization of the input capacity (i.e. cell supply function).  Further analysis of the TRM is proposed in Appendices~\ref{sec:REA} and \ref{sec:meq} where we show how the TRM can be interpreted as a REA algorithm and we analyze how it approximates shock waves (modified equations analysis).

In the CTM, the road is discretized into homogeneous cells of size $\Delta x$ (indexed in the ascending order in the direction of circulation of the road). The number of vehicles inside a cell $i$ is represented by a time-dependent function $\eta_i$. Consider some time step $\Delta t$. The evolution of the number of vehicles $\eta_i$ from a time $t$ to a time $t+\Delta t$  is described by the recurrence relation
\begin{equation}
 \eta_i(t+\Delta t)=\eta_i(t) + y_i(t) -  y_{i+1}(t)
 \label{eq:rec_ctm},
\end{equation}
where  the quantity $y_i(t)$ represents the number of vehicles entering the cell $i$ (from the cell $i-1$) between $t$ and $t+\Delta t$. These flows of vehicles are then assumed to be given, for any time $t$, by the relation
\begin{equation}
	y_i(t)=\min\lbrace \eta_{i-1}(t), Q_i(t), N_i(t)- \eta_{i}(t) \rbrace,
	\label{eq:flow_ctm}
\end{equation}
where $N_i(t)$ denotes the maximal number of vehicles the cell $i$ can hold at time $t$, and $Q_i(t)$ is the so-called input capacity of the cell $i$, which describes the maximal number of vehicles that could flow into the cell. Note in particular that for any $i$ and any $t$, $\eta_{i}(t)\in [0, N_i(t)]$.

Consider now the case where the maximal number of vehicles a cell can contain is constant, i.e. there exists some constant $N$ such that for any $i$ and any $t$, $N_i(t)=N$. Let us then take the input capacities $Q_i$ defined by
\begin{equation}
	Q_i(t)
	=\Delta t \; g\big(\frac{\eta_{i-1}(t)}{\Delta x},\, \frac{N-\eta_{i}(t)}{\Delta x}\big).
	\label{eq:input_cap}
\end{equation} 
Then, assuming that the CFL condition~\eqref{eq:CFL} is satisfied, we have $y_i(t)=Q_i(t)$ for any $i$ and any $t$.
Indeed, using the properties of $g$, we have on the one hand,
\begin{equation}
\begin{aligned}
	Q_i(t) 
	&= \Delta t \big( g\big(\frac{\eta_{i-1}(t)}{\Delta x},\, \frac{N-\eta_{i}(t)}{\Delta x}\big)
	-g\big(0,\, \frac{N-\eta_{i}(t)}{\Delta x}\big) \big) 
	 \le \Delta t K_1 \big\vert \frac{\eta_{i-1}(t)}{\Delta x} -0 \big\vert 
	= K_1\frac{\Delta t}{\Delta x}\eta_{i-1}(t),
\end{aligned}
\end{equation}
and similarly on the other hand,
\begin{equation}
\begin{aligned}
Q_i(t) 
= \Delta t \big( g\big(\frac{N}{\Delta x},\, \frac{N}{\Delta x}\big)
-g\big(\frac{N}{\Delta x},\, 0\big) \big)  \le K_2\frac{\Delta t}{\Delta x}\big(N- \eta_{i}(t)\big).
\end{aligned}
\end{equation}
Using then the fact that the CFL condition~\eqref{eq:CFL} implies that $K_1{\Delta t}/{\Delta x}\le 1$ and $K_2{\Delta t}/{\Delta x}\le 1$, we retrieve that $y_i(t)=Q_i(t)$.

	Note then that the density of vehicles $\rho_i(t)$ inside the $i$-th cell is obtained by dividing the number of vehicles $\eta_i(t)$ by the cell size $\Delta x$. Similarly, the maximal density of vehicles is $\rho_{\text{max}}= N/\Delta x$. Hence, by dividing \eqref{eq:rec_ctm} by $\Delta x$, we obtain the relation:
\begin{equation}
\rho_i(t+\Delta t)=\rho_i(t) + \frac{y_i(t)}{\Delta x} -  \frac{y_{i+1}(t)}{\Delta x}
\label{eq:ctm_trm}
\end{equation}
where  for any $i$,
\begin{equation*}
\begin{aligned}
\frac{y_i(t)}{\Delta x}
=\frac{Q_i(t)}{\Delta x}=
\frac{\Delta t}{\Delta x}\; g\big(\rho_{i-1}(t),\, \rho_{\text{max}} -\rho_{i}(t)\big), 
% \\
% &\frac{y_{i+1}(t)}{\Delta x}
% =\frac{Q_{i+1}(t)}{\Delta x}=
% \frac{\Delta t}{\Delta x}\; g\big(\rho_{i}(t),\, \rho_m -\rho_{i+1}(t)\big).
\end{aligned}
\end{equation*}
Hence, \eqref{eq:ctm_trm} coincide exactly with the recurrence relation defining the TRM \eqref{eq:discrete_kinetic}. We then conclude that the CTM,  when defined with the input capacities \eqref{eq:input_cap}, fits into the TRM framework proposed in this work. This particular choice of input capacities simplifies greatly the expression \eqref{eq:flow_ctm} defining the flow of vehicles: the minimum disappears and  we are left with the numerical fluxes of the TRM. %This means in particular that the physical interpretation of the TRM numerical flux given carries on to the CTM when defined with the input capacities \eqref{eq:input_cap}.
In conclusion, the TRM elegantly incorporates in a compact (analytic) form the restrictions on the flux of vehicles allowed to travel from cell to the other. %which ensure that we do not overflow a cell or remove non-existent vehicles from a cell: there is no need to distinguish between different cases anymore. 
The resulting smoothness opens the door for the application of results on polynomial systems of ODEs for the analysis of the numerical properties of the TRM.

% flatex input end: [Content/ctm.tex]

%% Interpretation as CTM

% flatex input end: [Content/num_exp.tex]

%--------------------------------------------------------------------------------------
%--------------------------------------------------------------------------------------

\section{Possible extensions of the TRM}\label{sec:ext}

In this section we present two straightforward extensions of the TRM that allow to model traffic in more complex situations that those covered by the LWR model~\eqref{eq:scalar_consv_pde} on a unidirectional road.

\subsection{Accounting for changing driving conditions}
The TRM can be adapted to model roads with non-homogeneous driving conditions, may they be due to changes in the maximal vehicle density admissible along the road (due for instance to varying number of lanes) or more generally to external factors affecting the optimal flow of vehicles. Indeed, consider once again a road discretized into cells $i \in I$, and let $ \rho_{\max}^{(i)}$ denote the capacity of the $i$-th cell. We call \textit{Extended TRM} the system of ODEs defined by 
	\begin{equation}
		\label{eq:ext_trm}
		\dot{\rho}_i(t) = \frac{1}{\Delta x} \big(C_{i}(t)F_{i}(\rho_{i-1}, \rho_{i}) - C_{i+1}(t)F_{i+1}(\rho_{i}, \rho_{i+1})
		+ R_{i}(\rho_i, t) - S_{i}(\rho_i, t)\big), \quad i \in I,\quad t\ge 0
	\end{equation}
where for each $i\in I$, the numerical flux $F_{i}$ is given by
\begin{equation} \label{eq:F}
	F_{i}(u, v) = g(u, \rho_{\max}^{(i)} -v),
\end{equation}
the function $g$ is defined in the same way as in \Cref{sec:trm}, and the function $t \mapsto C_{i}(t)$ takes values in $(0, 1]$ can be seen as local capacity drop factor which scales down the \q{ideal} flow of vehicles  $F_{i}(\rho_{i-1}, \rho_{i})$ between cells $i-1$ and $i$, i.e. the flow of vehicles one should expect between these cells given their density states. By identification with the usual TRM, this factor can essentially be seen as a time-varying normalized free flow speed $v_{\max}$ at the interface between cells $i-1$ and $i$, or within the kinetic compartmental interpretation, as a time-varying normalized reaction rate coefficient between compartments $i-1$ and $i$. On the other hand, allowing the maximal capacity $\rho_{\max}^{(i)}$ to change across compartments allows for instance to account for changes in the number of lanes or for non-uniform discretizations of the road.

Following the same arguments as the ones used in \Cref{thm:genFVM}, it is clear that the Extended TRM will yield solutions that preserve non-negativity and capacity.  Besides, the extended TRM has been shown to accurately describe real-world traffic datasets in applications related to traffic state estimation \cite{pereira2022parameter} and short-term prediction \cite{pereira2022short}. This is manly due to the local capacity drop factors $C_i$ which allow to locally change the conditions at which traffic flows on the roads, hence allowing the TRM (which inherits from the seemingly simplistic LWR model) to recreate complex traffic patterns observed in the data. 

\begin{figure}
    \centering
    \begin{subfigure}[t]{0.32\textwidth}
        \includegraphics[width=\textwidth,page=3,trim={0 1.5cm 0 3cm},clip]{fig_all.pdf}
        \caption{Representation of the road and its discretization} \label{fig:road}
    \end{subfigure}
        \begin{subfigure}[t]{0.32\textwidth}
        \includegraphics[width=\textwidth]{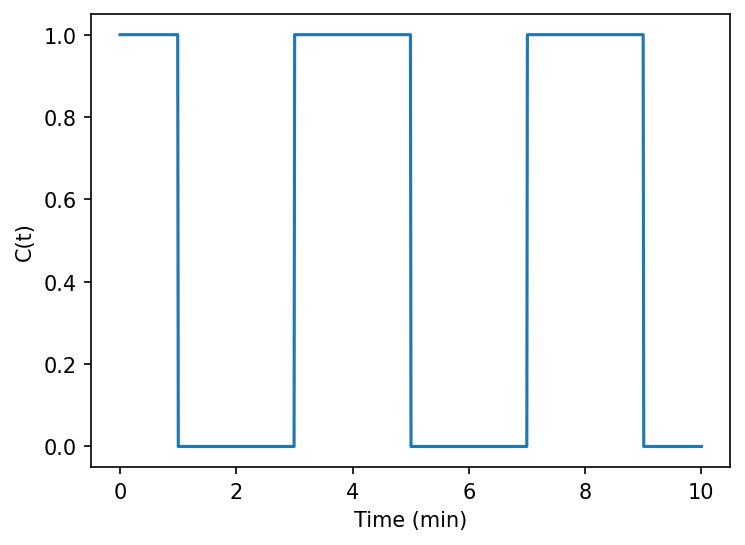}
        \caption{Time evolution of the capacity drop factor $C(t)$} \label{fig:cdrop}
    \end{subfigure} \\
    \begin{subfigure}{0.32\textwidth}
        \includegraphics[width=\textwidth]{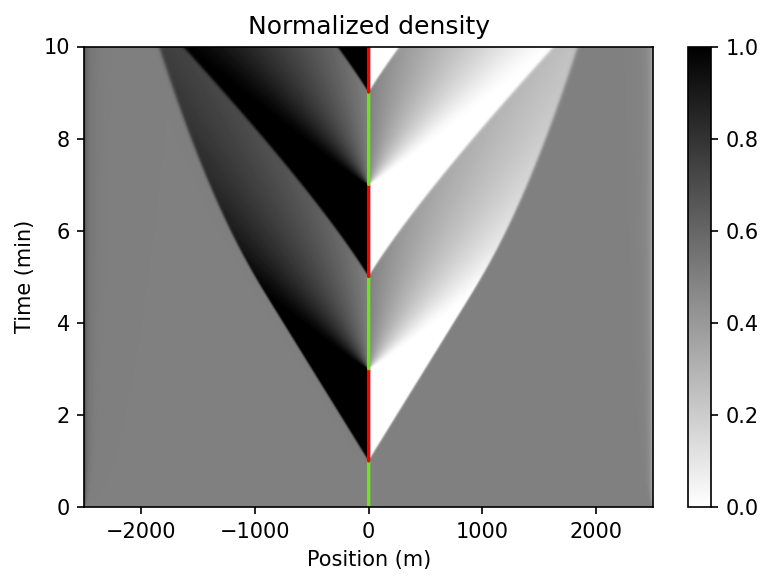}
        \caption{Space-time plot of the TRM density simulation} \label{fig:sim1}
    \end{subfigure}%
    \hfill%
    \begin{subfigure}{0.32\textwidth}
        \includegraphics[width=\textwidth]{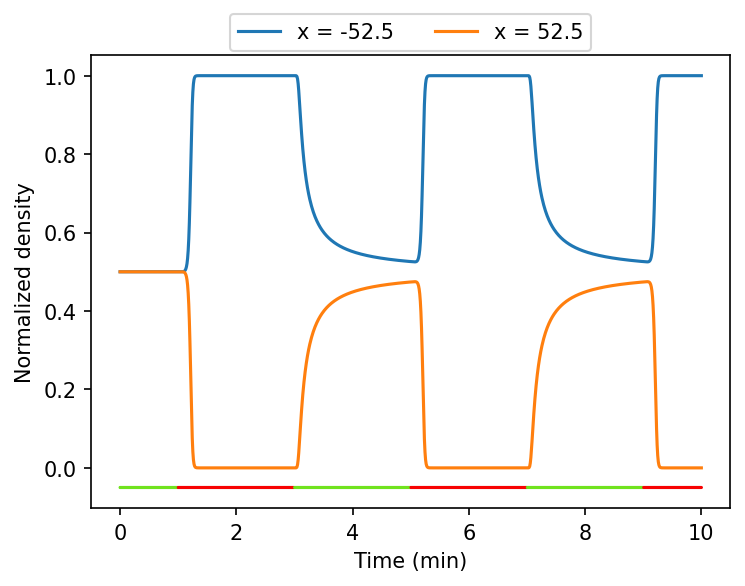}
        \caption{Time evolution of the density at $x=\pm 52.5$}\label{fig:sim2}
    \end{subfigure}%
    \hfill%
    \begin{subfigure}{0.32\textwidth}
        \includegraphics[width=\textwidth]{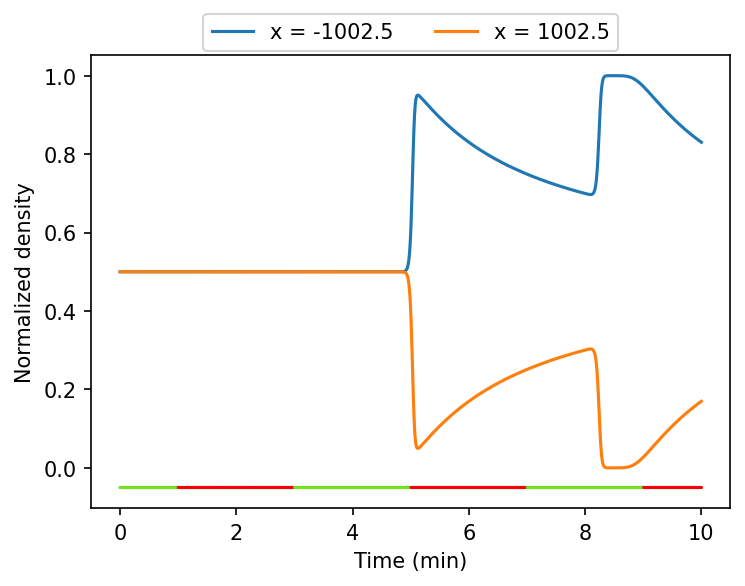}
        \caption{Time evolution of the density at $x=\pm 1002.5$}\label{fig:sim3}
    \end{subfigure}
    \caption{Extended TRM simulation of the normalized density on a road with a traffic light. The red and green lines indicate when the capacity drop factor is set to 0 (red lines) or to 1 (green lines).}
    \label{fig:res_sim}
\end{figure}

To illustrate this last property of the Extended TRM, we present in \Cref{fig:res_sim} the results of a simulation study where traffic along a unidirectional road with a traffic light is modeled. To do so, we start by discretizing the road into compartments  and consider that at the position $x=0$ where the traffic light is located, the corresponding capacity drop coefficient $C(t)$ oscillates between two values: 1 when the traffic light is green, and 0 when it is red (cf. \Cref{fig:road,fig:cdrop}). Following the kinetic interpretation of the TRM, this choice will effectively stop the flow of vehicles at $x=0$ when the traffic light is red. The mass action kinetic decomposition is chosen, which in particular yields a polynomial system of ODEs which can be solved at arbitrary time steps with high accuracy using standard libraries. 

The Extended TRM is run on a road of $5$km discretized into compartments of $5$m. The free flow speed of vehicles is set to $30$km/h and the simulation is run on a time horizon of 10min, with light switches happening approximately every 2 minutes. The results presented in \Cref{fig:sim1,fig:sim2,fig:sim3} seem to show that the model manages to recreate both the wave propagations and the oscillatory nature the density variations due to the traffic light changes.

\subsection{Accounting for road networks}

\begin{figure}
	\centering
	\includegraphics[height=0.2\textheight,page=1]{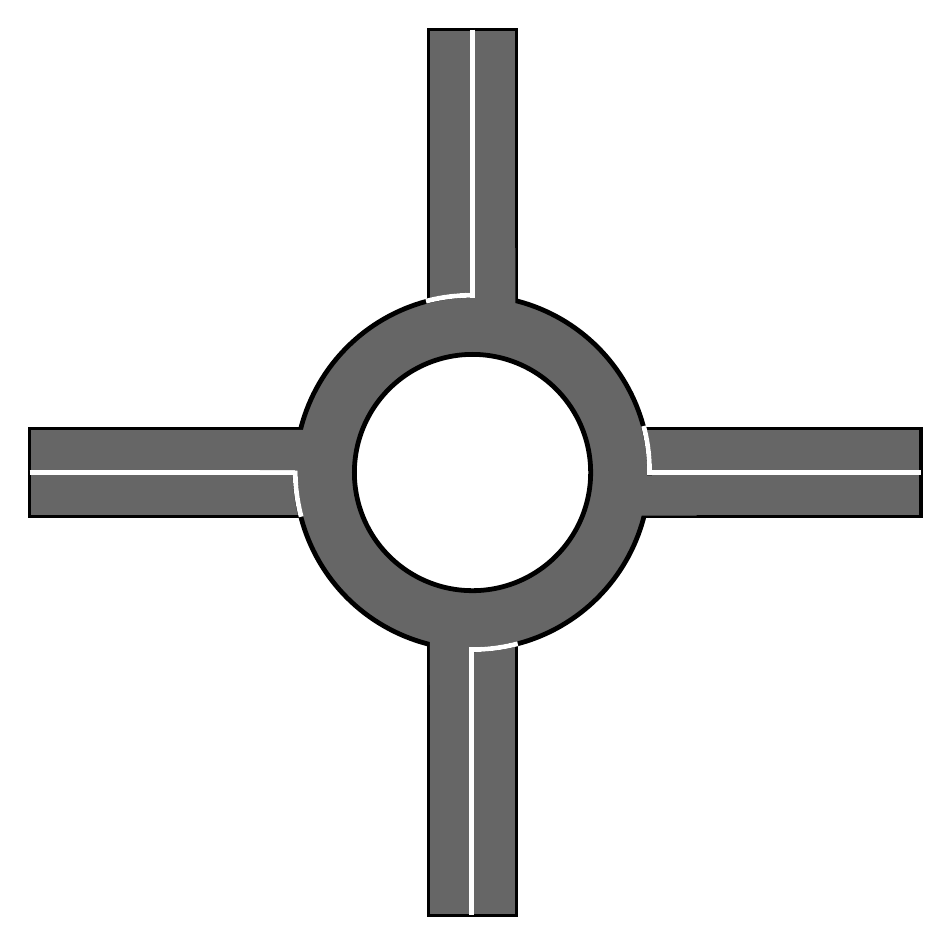}
	\hspace{1em}
	\includegraphics[height=0.2\textheight,page=2]{fig_all.pdf}
	\caption{TRM discretization \textit{(Right)} of a roundabout \textit{(Left)}. }\label{fig:rd_trm}
\end{figure}

The TRM can also be directly adapted to model traffic on networks. Indeed, following the kinetic and compartmental interpretation given in \Cref{sec:kinetic}, traffic on a network can be assimilated to a compartmental chemical reaction network for which each compartment is now allowed to have more than two other neighboring compartments. Starting from a road network, modeled as a directed graph, each directed edge is discretized into compartments, and each intersection is replaced by an additional compartment. An example of such discretization is given in \Cref{fig:rd_trm}, and show in particular how the TRM can be extended to handle complex intersections. The \textit{Network TRM} then takes the form of a system of ODEs defined as
	\begin{equation}
	\label{eq:net_trm}
	\dot{\rho}_i(t) = \frac{1}{\Delta x} \big(\sum_{j\in\mathcal N_{\mathrm{in}}(i)}F(\rho_{j}, \rho_{i}) - \sum_{k\in\mathcal N_{\mathrm{out}}(i)} F(\rho_{i}, \rho_{k})
	+ R_{i}(\rho_i, t) - S_{i}(\rho_i, t)\big), \quad i \in I,\quad t\ge 0
\end{equation}
where for each compartment $i\in I$, $\mathcal N_{\mathrm{in}}(i)$ denotes the set of compartments $j$ that point to $i$ (i.e. such that vehicles can move from $j$ to $i$) and $\mathcal N_{\mathrm{out}}(i)$ denotes the set of compartments $k$ that $i$ points to (i.e.  such that vehicles can move from $i$ to $k$).

\changednew{Once again, following the same arguments as the ones used in \Cref{thm:genFVM}, it is clear that the \emph{Network TRM} will also yield solutions that preserve non-negativity and capacity. This qualifies the Network TRM to act as a generic link and node model for traffic networks according to the classification of \cite{TAMPERE2011289}. Following the categories in \cite{TAMPERE2011289}, the \emph{Network TRM} extension is an unsignalized node model where supply and capacity constraints are handled in a balanced and continuous way. This latter property roots in fact in the balanced interaction of chemical species changing the concentrations in the compartments. In terms of supply constraint interaction rules \cite{TAMPERE2011289}, the \emph{Network TRM} generates them \emph{implicitly}: free space concentrations are transformed into occupied ones dynamically, smoothly. In terms of the solution algorithm, the set of nonlinear ODEs guarantees the existence of a unique solution to the density propagation accross all link and node compartments (via Lipschitz continuity \cite{Khalil}). Node and link models are described in a unified way with TRM. Links and nodes rely on the reaction rate governed exchange of material (vehicle) in- and outflows. Merging and diverging is handled by in- and outgoing flows (see, eq. \eqref{eq:net_trm}). 
Besides, the \emph{Network TRM} can readily be further extended in order to account for varying road conditions using the same approach as the one outlined for the Extended TRM. In practice, the link between the TRM and chemical reaction networks can be leveraged to deduce persistence and Lyapunov stability results that generalize even to time-varying networks, the properties of which are outlined in \Cref{sec:dyn} \citep{vaghy2022lyapunov,szederkenyi2022persistence,SzederVaghy_timevar}.}

\vspace{-0em}

\section{Conclusions}

Traffic Reaction Model (TRM), a family of Finite Volume Methods for segmenting certain hyperbolic Partial Differential Equations arising in traffic flow modelling, has been presented in the paper. First, the proposed numerical approximation scheme is consistent, monotone, nonnegative, capacitated, and conservative, even under non-zero sink and source terms. Second, the resulting system of nonlinear ODE models enables us to formally view traffic flow models as chemical reaction networks. %\changed{This new aspect allows us to apply the tools of compartmental systems and reaction network theory in the dynamical analysis of traffic models. Using this theory, a parameter-free logarithmic Lyapunov function was proposed for showing the robust stability of the dynamics of the ring topology.}
%Therefore, exactly the same nonnegative and kinetic ODE model can be obtained from two entirely different considerations. %This embeds TRM into a well established model class in system theory, which is an important advantage of the approach. It has been shown that the basic 'tube' model can be extended by on and off-ramps in a straightforward way still preserving the compartmental property. 
In the kinetic model formulation, two dual state variables appear representing traffic density and the density of free space, respectively. These dual densities enable a proper factorization of the nonlinearity in the original PDE, as done by the function $g()$, and also allow us to define multiple meaningful numerical fluxes.

We have compared the TRM to other known popular discretization schemes, and clarified their relations. Here we emphasize that TRM is equivalent to the CTM for a particular choice of input capacities. Numerical results show that the convergence properties of the proposed model are comparable to other popular schemes. The main benefits of the proposed modeling approach are the possibility to represent the system model in the form of smooth ODEs which is often required to apply advanced nonlinear control techniques, and to use the tools of compartmental systems and reaction network theory in the dynamical analysis of traffic models. The latter one was illustrated through the example of the ring topology, where the robust persistence and stability of the dynamics were shown using a Petri net representation (containing also both sets of state variables), and a parameter-free logarithmic Lyapunov function, respectively. %The persistence and stability analysis of the ring topology suggests that the results of chemical reaction network theory can be used in the future to support the dynamical analysis and synthesis of the traffic network models. 
Finally, we have  presented the possible extensions of the model firstly in the form of a time-varying nonlinear system, with basically the same physical interpretation, to describe capacity changes due to e.g., traffic lights, and secondly, as a straightforward network formulation which can be used to model an arbitrary road (directed graph) structure, thus allowing to better model real traffic data \cite{pereira2022parameter,pereira2022short}.

Further work will be focused on structured stability analysis, state estimation, and robust control design. \changednew{An interesting research direction can be to further analyze the dual variable triggered parametrization of $f(\rho)$ within the framework of variational theory in order to see how the Hamilton-Jacobi representation (see, \cite{LAVAL201317}) is changed.} 

%Additional future research directions encompass the extensions of the model (e.g. capacity drop), structured stability analysis, traffic controller, and observer design. 

% flatex input end: [Content/ccl.tex]

%% Empirical study

%--------------------------------------------------------------------------------------

\bibliographystyle{apalike} 
\bibliography{cleanbib}

%--------------------------------------------------------------------------------------

%%---------------------------------------------------------------------------

%--------------------------------------------------------------------------------------
%--------------------------------------------------------------------------------------

\begin{appendices}

\section{Well-posedness of the traffic flow PDE}\label{app:3}

\begin{lemma}\label{lem:unique}

Let $T>0$. Let $f:\Rone \rightarrow \Rone$ be bounded and locally Lipschitz-continuous, and let $\psi : [0,T]\times\Rone\times\Rone \rightarrow \Rone$ such that
\begin{itemize}
    \item $[(t,x) \in [0,T]\times\Rone \mapsto \psi(t,x,u)]$ is bounded and Lipschitz-continuous, uniformly in $u\in\Rone$,
    \item $[u\in\Rone \mapsto \psi(t,x,u)]$ is bounded and locally Lipschitz-continuous, uniformly in $(t,x)\in [0,T]\times\Rone$,
    \item the estimate $\vert \psi(t,x,u) \vert \le C(1+\vert u\vert)$ holds when $\vert u\vert$ is sufficiently large, for some $C>0$ independent of $t,x,u$.
\end{itemize}
Finally, let  $u_0\in L^1(\Rone)$ such that $u_0$ is bounded. Then, the Cauchy problem defined by
\begin{equation}\label{eq:Cauchy_gen}
\begin{cases}
    u_t +f(u)_x = \psi(t,x,u), & (t,x)\in]0,T]\times\Rone, \\
    u(0, x)=u_0(x), & x\in\Rone,
\end{cases}
\end{equation}
admits weak solutions, among which there is exactly one solution $u\in C([0,T]; L^1(\Rone))$ satisfying the so-called entropy inequalities: for any convex function $\eta\in C^2(\Rone)$ and function $q$ such that $\dot q=\dot \eta\dot f$, 
\begin{equation}
 \int_0^T\int_{\Rone}\big( \eta(u)\phi_t+q(u)\phi_x
-\dot\eta(u)\psi(t,x,u)\phi\big) dxdt
\ge 0
\label{eq:vv_sol_ineq}
\end{equation}
for any $\phi\in C^1_c([0,T]\times\Rone)$ such that $\phi\ge 0$. This particular solution is called entropy solution of the Cauchy problem \eqref{eq:Cauchy_gen}.
Moreover, the entropy solution $u$ depends continuously on the initial condition $u_0$ and is bounded.
\end{lemma}

\begin{proof}
This result a direct application of Theorems 3.1 and 6.1 of \cite{chen4quasilinear} in the case where the diffusion coefficient of the quasilinear anisotropic degenerate parabolic equation is zero.
\end{proof}

Note that since the entropy solution satisfies $u\in C([0,T];L^1(\Rone))$, it implies  that  for all $t_0 >0$,
$\lim\limits_{t\to t_0}\int_{\mathbb{R}}\vert u(
t,r)-u(t_0,r)\vert~dr=0,$ and $
 \lim\limits_{t\to 0+}\int_{\mathbb{R}}|u(t,r)-u_0(r)|\, dr=0.
$
Hence, the regularity in time (measured in the above $L_1$-norm sense) of the conserved quantity under nonzero sink and source terms is guaranteed.

% flatex input end: [Content/Appendix.tex]

\section{Empirical accuracy and convergence test for the TRM}\label{sec:numTRM}

Consider a source and sink term free traffic flow model. The empirical accuracy tests are carried out with the parameters $\omega=1$, $\rho_{\text{max}}=100$, and $\Delta x = L / P$, where $L=20$ is the length of the examined segment, and $N$ is the number of discrete segments. The end time is $T=2/60$. We also consider that we have reflexive boundary conditions, i.e. we set $\rho_{-1} = \rho_0$ and $\rho_{P+1} = \rho_P$. We compare the performances of the semi-discrete scheme \eqref{eq:gen_scheme} and of the fully-discrete scheme \eqref{eq:discrete_kinetic} in terms of accuracy, for increasing numbers of space discretization points $P$. First, the system of ODEs of the semi-discrete scheme is solved numerically (Runge-Kutta approximation). Second, the time step of the fully discrete scheme is  systematically set to $\Delta x / (2v_{\text{max}})$ (following \Cref{thm:genFVM2}).

First, we introduce the time dependent error function $e$, which is the spatial $L^1$-norm error in each time $t$ as follows
\begin{equation}
\label{eq:error}
e(t) = \sum_{i=1}^{P}\int_{x_{i-1/2}}^{x_{i+1/2}}\left\vert\rho(x,t) - \rho_i(t) \right\vert dx,
\end{equation}
where $\rho(x,t)$ denotes the value of the true solution of the PDE and $\rho_i(t)$ denotes, for the $i$-th cell, either the solution of the ODE system of the semi-discrete scheme or a piecewise constant reconstruction obtained from the quantities computed by the fully discrete scheme (and equal to $\rho_i^n$ in each interval $[t_n, t_{n+1}]$). This spatial error term is computed analytically for each $t$ and arbitrary values of $\rho_i(t)$, in the case of Riemann problems, using the closed forms for $\rho(x, t)$ provided by \cite{Book_Leveque1992}.
After that, we can define the $L^1$- and $L^\infty$-norm of the spatial error function $e$ as follows
\begin{equation*}
\norm{e}_{1} = \int_0^{T} e(t)\, dt~~ \text{and}~~ \norm{e}_{\infty} = \max_{0\leq t  \leq T} e(t).
\end{equation*}

%\BK{Nekem nem vilagos a terbeli es idobeli norma kulonbsege. Be tudjuk mind a kettot vezetni?}
%Both quantities are computed using Matlab numerical integrators and optimizers. % and calculate the norms in Eqs. \eqref{eq:norm1}-\eqref{eq:norminf}.
\subsection{Shock wave}
\begin{figure}
	\centering
		\includegraphics[width=0.6\textwidth]{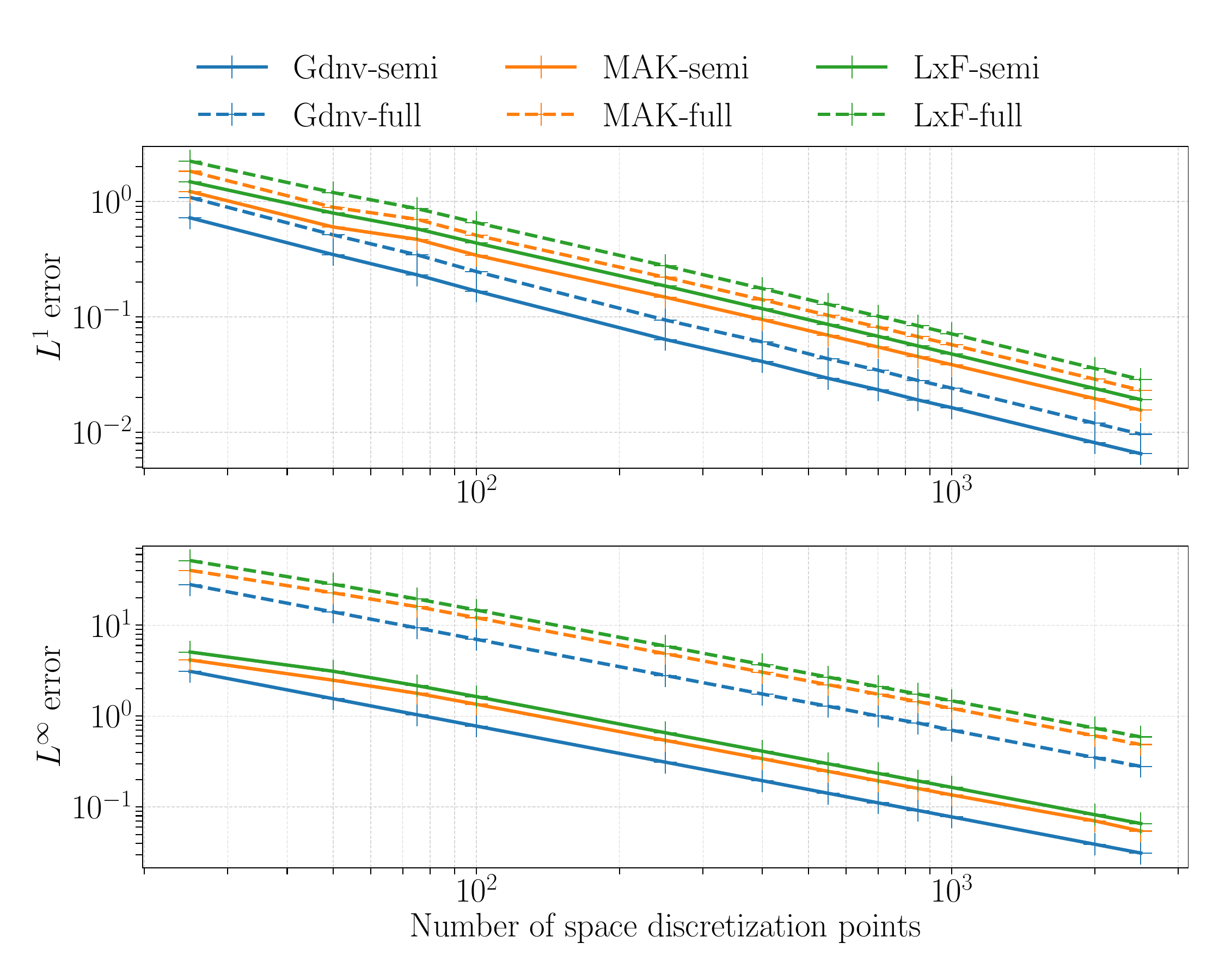}
  \caption{$L^1$-norm and $L^\infty$-norm errors  of the shock wave discretization for various numbers of segments  $P$.}\label{fig:shock}
\end{figure}

Consider the Riemann-problem with the initial condition
\begin{equation*}
\rho(x, 0) =     
\begin{cases}
0.1\rho_{\text{max}}  & \text{if } x < L/2 \\
0.8\rho_{\text{max}} & \text{otherwise},
\end{cases}
\end{equation*}
which induces a solution consisting of a shock wave.
Figure \ref{fig:shock} shows the discretization errors of various monotone discretization schemes (both semi and fully discrete): the mass action kinetic ($\Mk$) TRM \eqref{eq:g1}, the Godunov (Gdnv) scheme \eqref{eq:g2}, and the modified Lax-Friedrichs (LxF) scheme 
 \footnote{Numerical flux given by $F(u, v) = (f(u) + f(v))/2 + d(u - v)$ for $d \geq \omega \rho_{\text{max}} / 2$. Note that LxF is not a kinetic scheme.}. The errors decrease approximately as $\mathcal{O}(P^{-1})$ for all (semi and fully discrete) schemes, and the semi-discrete formulation of a given scheme systematically over-performs its fully discrete formulation. This is expected as the fully-discrete formulation add additional error due to the time discretization. If we compare the semi or fully discrete schemes with one another, the Godunov scheme systematically overperforms the other two, which display similar errors.

\subsection{Rarefaction wave}

\begin{figure}
	\centering
		\includegraphics[width=0.6\textwidth]{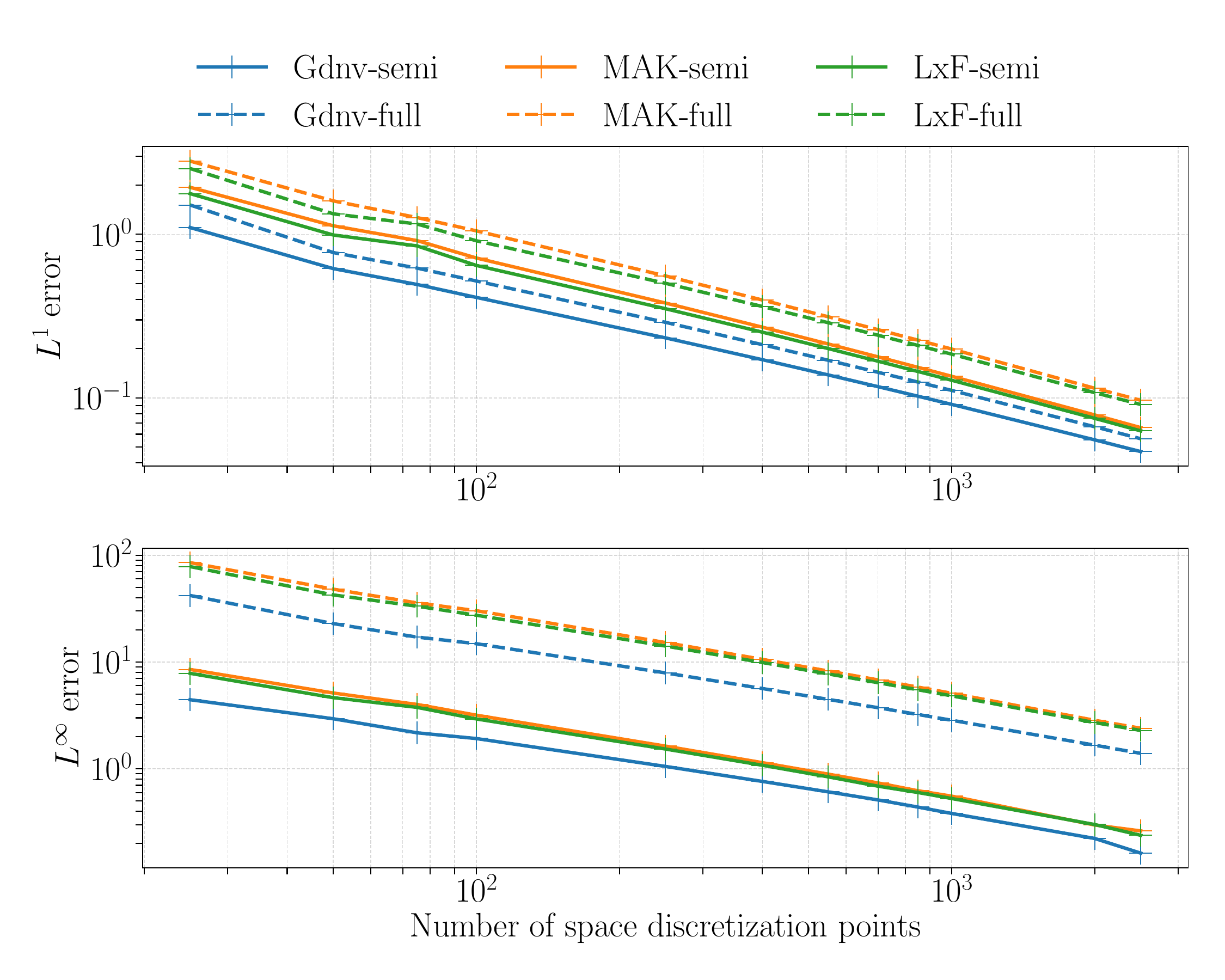}
  \caption{$L^1$-norm and $L^\infty$-norm errors  of the rarefaction wave discretization for various numbers of segments  $P$.}\label{fig:rar}
\end{figure}

Consider the Riemann-problem with the initial condition
\begin{equation*}
\rho(x, 0) = 
\begin{cases}
0.8\rho_{\text{max}} & \text{if } x < L/2 \\
0.1\rho_{\text{max}} & \text{otherwise},
\end{cases}
\end{equation*}
which induces this time a solution consisting of a rarefaction wave.
The discretization errors associated with this choice of initial condition are given in Figure \ref{fig:rar}. Once again, the semi-discrete formulations overperform compare to the fully-discrete ones, and the Godunov scheme overperforms compared to the other two schemes, that display similar errors. But this time, the errors at a slightly lower order of $\mathcal{O}(P^{-3/4})$. This difference could be explained by the fact that (spatially) piecewise constant reconstructions are used here to approximate a continuous rarefaction fan.

\section{Interpretation of the TRM as a REA algorithm}\label{sec:REA}

Let us consider the fully discrete TRM with no source terms. Then, if $\rho_i^k$ denotes the density of vehicles  and  $\nu_i^k=\rho_{\max}-\rho_i^k$ denotes the density of empty space in the $i$-th cell at time $t_k$ (for any $i\in\mathbb{Z}$, $k\in\mathbb{N}$), we have 
\begin{equation*}
	\begin{aligned}
		\rho_i^{k+1}&=\rho_i^{k}+\frac{\Delta t}{\Delta x}\left(g(\rho_{i-1}^{k},\nu_{i}^{k})-g(\rho_{i}^{k},\nu_{i+1}^{k})\right) \\
		&=\rho_i^{k}+\frac{\Delta t}{\Delta x}\left(g(\rho_{i-1}^{k},\nu_{i}^{k})-g(\rho_{i}^{k},\nu_{i}^{k})\right)
		+\frac{\Delta t}{\Delta x}\left(g(\rho_{i}^{k},\nu_{i}^{k})-g(\rho_{i}^{k},\nu_{i+1}^{k})\right),
	\end{aligned}
\end{equation*}
Let us then define the quantities $\alpha_l(i)$ and $\alpha_r(i)$ by
\begin{equation}\label{eq:advsys_coiefs}
	\begin{aligned}
		\alpha_l(i)
		=\frac{g(\rho_{i-1}^{k},\nu_{i}^{k})-g(\rho_{i}^{k},\nu_{i}^{k})}{\rho_{i-1}^k - \rho_i^k}, \quad  
		\alpha_r(i)
		=\frac{g(\rho_{i}^{k},\nu_{i}^{k})-g(\rho_{i}^{k},\nu_{i+1}^{k})}{\nu_{i+1}^k - \nu_{i}^k}, 
	\end{aligned}
\end{equation}
where by convention we take $\alpha_l(i)=0$ (resp. $\alpha_r(i)=0$) if $\rho_{i-1}^k = \rho_i^k$ (resp. $\nu_{i+1}^k = \nu_{i}^k$). Note in particular that since $g$ is non-decreasing with respect to each of its argument, we have $\alpha_l(i)\ge 0$ and $\alpha_r(i)\le 0$.
We can then write the TRM recurrence relation as
\begin{equation}\label{eq:update_rho}
	\begin{aligned}
		\rho_i^{k+1}=\frac{1}{2}\big[  I^+\left(\alpha_l(i),\rho_{i-1}^k,\rho_i^k\right)
		+\left(\rho_{\text{max}}-I^-\left(\alpha_r(i),\nu_{i}^k,\nu_{i+1}^k\right)\right)
		\big],
	\end{aligned}
\end{equation}
where
\begin{equation}\label{eq:advsys}
	\begin{aligned}
		I^+\left(\alpha_l(i),\rho_{i-1}^k,\rho_i^k\right) 
		&=\rho_i^{k}-\frac{\Delta t}{\Delta x/2}\alpha_l(i)(\rho_{i}^k-\rho_{i-1}^k),
		\\
		I^-\left(\alpha_r(i),\nu_{i}^k,\nu_{i+1}^k\right)
		&=\nu_i^{k}-\frac{\Delta t}{\Delta x/2}\alpha_r(i)(\nu_{i+1}^k-\nu_{i}^k). 
	\end{aligned}
\end{equation}
Hence, the update of the vehicle density within the $i$-th cell can be written as the average between two quantities: one depending only on the vehicle density of the cells $i-1$ and $i$, and one depending only on the density of empty space of the cells $i$ and $i+1$. It turns out that these quantities can be linked to the solutions of some specific advection problems.

Indeed, consider the $i$-th cell, centered at $x-i$, its left border $x_{i-1/2}$, and right border $x_{i+1/2}$ (cf. \Cref{fig:sk_half}). Let us then introduce the functions $\eta_l$ and $\eta_r$ respectively centered around the left and right border of the cell, and given by
\begin{equation*}
	\eta_l(x)=\begin{cases}
		\rho_{i-1}^k & \text{if } x<x_{i-1/2} \\
		\rho_{i}^k & \text{if } x>x_{i-1/2}
	\end{cases}, 
	\quad
	\eta_r(x)=\begin{cases}
		\nu_i^k & \text{if } x<x_{i+1/2} \\
		\nu_{i+1}^k & \text{if } x>x_{i+1/2}
	\end{cases}.
\end{equation*}
Consider the advection problem defined by
\begin{equation*}
	\left\lbrace\begin{aligned}
		&\rho_t + \alpha_l(i)\rho_x =0, \quad x \in \R, \quad t\ge t_k, \\
		&\rho(0,x)=\eta_l(x), \quad x \in\R.
	\end{aligned}\right.
\end{equation*}
Then, the average over the \textit{left half} of the $i$-th cell of the solution at time $t_{k+1}=t_k+\Delta t$ of this advection problem is given by the quantity $I^+\left(\alpha_l(i),\rho_{i-1}^k,\rho_i^k\right)$ (cf. Appendix \ref{sec:adv}). 
Similarly, the quantity $I^-\left(\alpha_r(i),\nu_{i}^k,\nu_{i+1}^k\right)$ corresponds exactly to the cell average over the \textit{right half} of the $i$-th cell
of the solution at time $t_{k+1}$ of the advection problem with coefficient $\alpha_r(i)\le 0$ and (piecewise constant) initial condition $\eta_r$.

\begin{figure}
	\centering
	\resizebox{0.7\textwidth}{!}{
		% flatex input: [sk_half.tex]

		\begin{tikzpicture}[x=\textwidth,y=0.25\textwidth,every node/.style={scale=1.35}]

			\draw[thick] (0.01,0.75) -- (0.99,0.75);
			\draw[thick](0.01,0.45) -- (0.99,0.45);

			\draw[dashed] (0.05,1) -- (0.05,0.1);
			\draw[dashed] (0.35,1) -- (0.35,0.1);
			\draw[dashed] (0.65,1) -- (0.65,0.1);
			\draw[dashed] (0.95,1) -- (0.95,0.1);
			
			\draw[] (0.05,0.75) -- (0.05,0.45);
			\draw[dashed] (0.35,0.75) -- (0.35,0.45);
			\draw[dashed] (0.65,0.75) -- (0.65,0.45);
			\draw[dashed] (0.95,0.75) -- (0.95,0.45);
			
			\draw (0.2,0.9) node {$(i-1)$-th cell};
			\draw (0.5,0.9) node {$i$-th cell};
			\draw (0.8,0.9) node {$(i+1)$-th cell};
			
			\draw[-{Latex[length=3mm,width=2mm]}] (0,0.15) -- (1,0.15) node[anchor=west]{$x$};
			\draw (0.2,0.12) node[anchor=north]{$x_{i-1}$} -- (0.2,0.18) ;
			\draw (0.5,0.12) node[anchor=north] (tot) {$x_{i}$} -- (0.5,0.18) ;
			\draw (0.8,0.12) node[anchor=north]{$x_{i+1}$} -- (0.8,0.18) ;
			
			\draw (0.35,0.12) node[anchor=north]{$x_{i}-\Delta x/2$} -- (0.35,0.18) ;
			\draw (0.65,0.12) node[anchor=north]{$x_{i}+\Delta x/2$} -- (0.65,0.18) ;

			\draw (0.425,0.6) node {Left half};
			\draw (0.575,0.6) node {Right half};

			\draw[dashed,red] (0.5,0.15) -- (0.5,0.75);
			\draw[dashed,red] (0.8,0.15) -- (0.8,0.75);
			\draw[dashed,red] (0.2,0.15) -- (0.2,0.75);

			%\draw[decorate,decoration={brace,mirror}] (0.35,0) -- (0.5,0) node[below,midway,yshift=-5pt] { $J_l(t_n,x_i)$};
			
			%\draw[decorate,decoration={brace,mirror}] (0.5,-0) -- (0.65,-0) node[below,midway,yshift=-5pt] { $J_r(t_n,x_i)$};

		\end{tikzpicture}

		% flatex input end: [sk_half.tex]
		
	}
	\caption{Left and right halves of the cells.}
	\label{fig:sk_half}
\end{figure}

Hence, the TRM can be seen as a particular instance of Reconstruct-Evolve-Average (REA) algorithm \cite{leveque2002finite}. Indeed, assuming the density values at time $t_k$ are known, the values at time $t_{k+1}$ can be seen as resulting from the following three-step process:
\begin{description}
	\item[Step R] Reconstruct a density function $\rho(t_k, \cdot)$ which is constant over each cell from the current cell estimates $\rho_i^k$ and define the free space as $\nu(t_k, \cdot)=\rho_{\text{max}}-\rho(t_k, \cdot)$
	\item[Step E] Evolve, at each inter-cell boundary $x_{i+1/2}$, the system of (decoupled) advection equations given by
	\begin{equation}\label{eq:adv_syst}
		\left\lbrace
		\begin{aligned}
			&\rho_t + \alpha_l(i+1)\rho_x =0 \\
			&\nu_t + \alpha_r(i)\nu_x=0
		\end{aligned}
		\right.
	\end{equation}
	from the initial conditions $(\rho(t_k, \cdot), \nu(t_k, \cdot))$, and for a time $\Delta t$.
	\item[Step A] Compute the cell estimate $\rho_i^{k+1}$ of the $i$-th cell by averaging the value of the cell average of $\rho$ over the left half of this cell, and the value of  cell average of $(\rho_{\text{max}} - \nu)$ over its right half.
\end{description}

REA algorithms include the Godunov scheme \cite{leveque2002finite}, which can be described as an implementation of the update rule \eqref{eq:update_rho} where the advection problems are replaced by coupled system of conservation laws 
\begin{equation}\label{eq:adv_syst1}
	\left\lbrace
	\begin{aligned}
		&\rho_t + g(\rho,\nu)_x =0 \\
		&\nu_t -g(\rho,\nu)_x=0.
	\end{aligned}
	\right. \,
\end{equation}
This adjoint system can actually be decoupled and then yields that $\rho$ satisfies our initial PDE (without no sink and source terms)
\begin{equation*}\label{eq:pde}
	\partial_t\rho+\partial_xf(\rho)_x =0,
\end{equation*}
and the relation $\nu=\rho_{\text{max}}-\rho$. Hence, we actually proved that solving the resolution of the adjoint system
as implemented in the Godunov scheme is equivalent to that of advection problems with coefficients as defined in \eqref{eq:advsys} and 
\begin{equation*}
	g=g_\Gd(\rho, \nu) = \min(D(\rho), Q(\rho_{\max} - \nu)).
\end{equation*}

Other choices of decompositions $g(\rho,\nu)$ result in different advection coefficients and consequently may trigger different numerical approximation errors. A systematic and comparative analysis of different $g(\rho,\nu)$ functions can be obtained by evolving the time-space bounds to their respective advection coefficients and compare them.

\section{Error analysis through modified equations}\label{sec:meq}

In this section we show how the smoothness of the numerical fluxes of the TRM for the Greenshields flux and defined from the mass action kinetic decomposition $g(\rho, \nu)=g_\Mk(\rho, \nu) = \omega \rho \nu$ can be leveraged to provide a priori estimates on the approximation of shocks by the scheme. In particular, we give a formula modeling the smooth shocks obtained from the TRM and compare them to the one obtained using a different scheme (Lax--Friedrichs).

Let us consider the fully discrete TRM defined from the mass action kinetic decomposition, and assume that there exists some constant $\delta = (\Delta t/\Delta x)v_{\max} \in (0,1/2)$. Then, the recurrence relation defining this scheme can be written in the following form (for  $i\in\mathbb{Z}$,  $k\in\mathbb{N}_0$)
\begin{equation}\label{eq:trunc_1}
	\rho_i^{k+1}-\rho_i^{k}
	-\frac{\Delta t}{\Delta x}\big(F(\rho_{i-1}^{k},
	\rho_{i}^{k})-F(\rho_{i}^{k},\rho_{i+1}^{k})\big)=0.
\end{equation}

Replacing then,  in the left-hand side of~\eqref{eq:trunc_1}, the estimates $\rho_j^n$ by the corresponding evaluations $\rho(t_n, x_j)$ of the PDE solution $\rho$, yields the so-called local truncation error (LTE) $L_\rho(t,x)$, given for any $t\ge0, x\in\R$ by:
\begin{equation}\label{eq:lte_def}
	\begin{aligned}
		L_\rho(t,x)=~ &\rho(t+\Delta t, x) -\rho(t,x)\\
		&-\frac{\Delta t}{\Delta x}\big( F(\rho(t,x-\Delta x),\rho(t,x))
		-F(\rho(t,x),\rho(t,x+\Delta x))\big) .
	\end{aligned}
\end{equation}
The LTE is the difference between the true solution, one time step ahead, and the approximation of this quantity using the recurrence relation of the scheme. Hence, it provides an indication of how well the recurrence relation defining a numerical scheme approximates the variations (in time) of the true solution. Indeed, if $L_\rho(t,x)=0$ for any $t\ge 0, x\in\R$, and assuming that  $\rho_i^0=\rho(0, x_i)$ for any $i\in\mathbb{Z}$,  then the approximate solution TRM recurrence~\eqref{eq:trunc_1} will coincide with the true solution at any time $t_k$ and location $x_i$.

Building on this idea, the modified equation approach \cite{Book_Leveque1992} proposes to find smooth functions $\tilde{\rho}$, for which estimates of the magnitude of LTE (and hence of the ability of the scheme to approximate $\tilde{\rho}$) are available. This done by expanding the LTE $L_{\tilde{\rho}}$ in~\eqref{eq:lte_def} using a Taylor expansion of $\tilde{\rho}$, and truncating the result at a given order $n$. Taking then $\tilde{\rho}$ to be the solution of the PDE defined by the remaining terms of the expansion yields that the corresponding LTE $L_{\tilde{\rho}}$ must be of order $\mathcal{O}(\Delta x^n)$. For instance, when truncating at an order $1$, we get that $\tilde{\rho}=\rho$ is the solution of the conservation law~\eqref{eq:pde} (with no source terms), and therefore that the LTE of the TRM for $\rho$ satisfies $L_{\rho}(t,x)=\mathcal{O}(\Delta x)$ (which is in accordance with the fact that monotone schemes for conservation laws converge at most as fast as $\mathcal{O}(\Delta x)$ \cite{leveque2002finite}).

Considering a truncation of order $2$ yields the following PDE to define $\tilde{\rho}$,
\begin{equation}\label{eq:pde_mlte}
	\partial_t\tilde{\rho} + \omega \partial_x\big(\tilde{\rho}(\rho_{\max}-\tilde{\rho})\big) + \frac{\Delta t}{2}\partial_{tt}\tilde{\rho}-\frac{\Delta x}{2}\omega\rho_{\max}\partial_{xx}\tilde{\rho}=0, 
\end{equation}
for which we know that the corresponding LTE $L_{\tilde{\rho}}$ satisfies $L_{\tilde{\rho}}(t,x)=\mathcal{O}(\Delta x^2)$. It is then natural to consider, for simple cases, the solutions of the so-called modified equation (of order 2)~\eqref{eq:pde_mlte}, in order to get a better understanding of the numerical properties of the TRM scheme, since the scheme yields better approximations of the solutions of this PDE ($L_{\tilde{\rho}}(t,x)=\mathcal{O}(\Delta x^2)$) than of the solutions of the conservation law ($L_{\rho}(t,x)=\mathcal{O}(\Delta x)$). 

One such case is the study of the approximation of shocks by the numerical scheme. When the initial condition of the conservation law contains a discontinuity such that the left value $\rho_l$ is less than the right value $\rho_r$, then the solution consist of a shock traveling with a speed $v_s=v_{\max}(1-(\rho_l+\rho_r)/\rho_{\max})$. A similar behavior is observed for some solutions $\tilde{\rho}$ of~\eqref{eq:pde_mlte}, as stated in the next proposition (proof in Appendix \ref{sec:proof_meq}).

\begin{proposition}\label{prop:modeq_shock}
	Let $u$ denote the sigmoid function, which is defined by $u(z)=(1+\exp(-z))^{-1}$, $z\in\R$.
	Let $\rho_l, \rho_r \in [0, \rho_{\max}]$ such that $\rho_l < \rho_r$ and let $\tilde\rho : \R_+\times\R \rightarrow \R$ be defined by
	\begin{equation}\label{eq:tilderho}
		\tilde{\rho}(t,x)=\rho_l + (\rho_r-\rho_l)u\big(\sigma^{-1}(x-x_0-v_st)\big)
	\end{equation}
	where $x_0\in\R$ is an arbitrary constant and
	\begin{equation*}
		\sigma = \sigma_{\text{TRM}}=\frac{\rho_{\max}}{2(\rho_r-\rho_l)} \bigg(1-\delta\bigg(1-\frac{\rho_l+\rho_r}{\rho_{\max}}\bigg)^2\bigg)\, \Delta x
	\end{equation*}
	with $\delta = ({\Delta x}/{\Delta t})v_{\max} \in (0, 1/2)$. 
	Then, $\tilde{\rho}$ is the only traveling wave solution of~\eqref{eq:pde_mlte} satisfying $\lim\limits_{x\rightarrow -\infty} \tilde\rho(0, x)=\rho_l$
	and $\lim\limits_{x\rightarrow +\infty} \tilde\rho(0, x)=\rho_r$.
\end{proposition}

Note that the traveling wave solution described in \Cref{prop:modeq_shock} travels at the same speed as the shocks of PDE~\eqref{eq:pde}. The wave in question consists in a smooth approximation of a shock from $\rho_l$ to $\rho_r$ using a sigmoid function. As for the arbitrary constant $x_0$, it sets the position, at time $t=0$, of the inflection point of the traveling wave. 
The parameter $\sigma$, on the other hand, sets how sharp the transition between $\rho_l$ and $\rho_r$ is: indeed, for any $\epsilon\in (0,1)$, the errors $\vert \rho_l - \tilde\rho\vert$ and $\vert \rho_r - \tilde\rho\vert$ are smaller than $\epsilon(\rho_r -\rho_l)$ at a distance (from the inflection point) greater than $\sigma \vert\log(\epsilon/(1-\epsilon))\vert$. For instance,  errors smaller than $0.01(\rho_r -\rho_l)$ are achieved $5\sigma$ away from the inflection point.

Consequently, the quantity $\sigma$ can be used to  get an idea of the extent with which the TRM smoothes out the discontinuities. In particular, since $\sigma$ is proportional to $\Delta x$, the expression of $\sigma$ given in \Cref{prop:modeq_shock} can be used to directly estimate the number of cells needed to transition from the left state of a shock to the right state. Also, since $\sigma$ is inversely proportional to the shock size $\rho_r-\rho_l$, we expect this transition to be longer when approximating small shocks than when approximating bigger shocks.

Note that similar conclusions can be drawn for the Lax--Friedrichs (LxF) scheme. Indeed, using the same approach, one can show that the modified equation of order $2$ of the LxF scheme takes the form
\begin{equation}\label{eq:pde_mltelf}
	\partial_t\tilde{\rho} + \omega \partial_x\big(\tilde{\rho}(\rho_{\max}-\tilde{\rho})\big) + \frac{\Delta t}{2}\partial_{tt}\tilde{\rho}-\frac{\Delta x^2}{2\Delta t}\omega\rho_{\max}\partial_{xx}\tilde{\rho}=0,
	% \tilde{\rho}_t + \omega \big(\tilde{\rho}(\rho_{\max}-\tilde{\rho})\big)_x + \frac{\Delta t}{2}\tilde{\rho}_{tt}-\frac{\Delta x^2}{2\Delta t}\tilde{\rho}_{xx}=0, 
\end{equation}
which differs from the TRM only by the factor in front of $\tilde\rho_{xx}$. Hence, the same reasoning as the one used in \Cref{prop:modeq_shock} can be applied, and we can deduce that the traveling wave solution of \eqref{eq:pde_mltelf} can once again be written as in ~\eqref{eq:tilderho}, but with a parameter $\sigma=\sigma_{\text{LxF}}$ now equal to
\begin{equation}
	\sigma_{\text{LxF}}=\frac{\rho_{\max}}{2(\rho_r-\rho_l)} \bigg(\frac{1}{\delta}-\delta\bigg(1-\frac{\rho_l+\rho_r}{\rho_{\max}}\bigg)^2\bigg)\, \Delta x.
\end{equation}
In particular, since $\delta \in (0,1/2)$ we have that $\sigma_{\text{TRM}} < \sigma_{\text{LxF}}$, and therefore, we can expect the mass action kinetic TRM to smooth less the shocks than the LxF scheme.

\begin{figure}[h]
	\centering
	\includegraphics[width=0.6\textwidth]{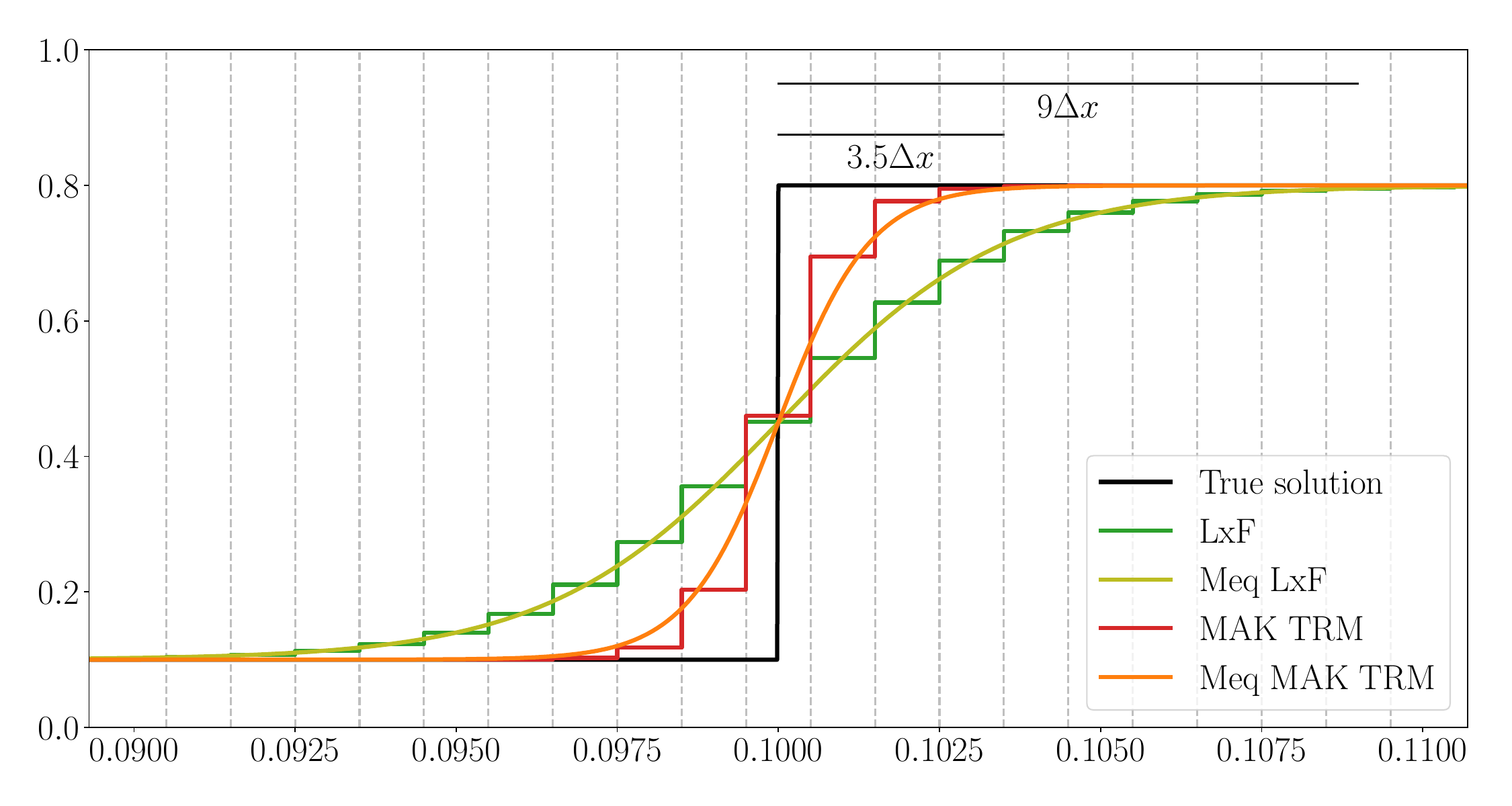}
	\caption{Comparison of the Lax--Friedrichs (LxF) scheme and the mass action kinetic TRM. The true solution (in black) corresponds to PDE~\eqref{eq:pde} with no source and sink term, and initial condition equal to $\rho_l=0.1$ if $x<0$ and $\rho_r=0.9$ if $x>0$ ($\rho_{\max}=1, v_{\max}=1$). The solutions are computed at time $t=1$, and with parameters $\Delta x=10^{-3}$ and $\delta = (\Delta x / \Delta t)v_{\max}=0.4$. The LxF solution is denoted by LxF, the solution of modified equation of LxF is denoted by Meq LxF, the mass action kinetic TRM solution is denoted by MAK TRM and the solution of the modified equation of the mass action kinetic TRM is denoted by Meq MAK TRM. }\label{fig:meq0}
\end{figure}

Finally, we present in \Cref{fig:meq0} a comparison between approximations of a solution of PDE~\eqref{eq:pde} using the mass action kinetic TRM  and the LxF scheme. We also plot the solutions of the corresponding modified equations of order 2. As expected, the mass action kinetic TRM and LxF approximate solutions are close to their respective solutions to the modified equation. Since we consider $\rho_l=0.2$ and $\rho_r=\rho_m$, we expect the mass action kinetic TRM solution to be close to the left and right state if we are at a distance of more than $5\sigma_{\text{TRM}}\approx 3.5\Delta x$ from the inflection point. As for the LxF scheme, we expect this distance to be $5\sigma_{\text{LxF}}\approx 9\Delta x$. The figure confirms these estimates.

% flatex input end: [Content/meq.tex]	

	%%---------------------------------------------------------------------------
	
	\section{Solutions of advection problems}\label{sec:adv}
	% flatex input: [Content/appen_adv.tex]

	Let $w_l, w_r \in\R$, and define $w_0$ by
	\begin{equation}
		w_0(x)=\begin{cases}
			w_l & \text{if } x<0 \\
			w_r & \text{if } x>0
		\end{cases}
		\label{ic_pc}
	\end{equation}
	For $\alpha\in\R$, consider the following advection problem:
	\begin{equation}
		\left\lbrace 
		\begin{aligned}
			& w_t +\alpha w_x=0, \quad t\in\R_+, \quad x\in\R \\
			& w(0,\cdot)=w_0
		\end{aligned}\right.
		\label{adv}
	\end{equation}
	The solution of \eqref{adv} is given by
	\begin{equation*}
		w(t,x)=w_0(x-\alpha t), \quad t\in\R_+, \quad x\in\R,
	\end{equation*}
	and therefore consists of a single shock (with left state $w_l$ and right state $w_r$) traveling at speed $\alpha$.
	
	Let then  $\Delta t, \Delta x >0$ and let us now assume that $\alpha$ satisfies
	\begin{equation}
		\vert \alpha \vert \Delta t \le \frac{\Delta x}{2}
		\label{ccfl}
	\end{equation}
	Let us then introduce the quantities $I^+(\alpha,w_l,w_r)$ and $I^-(w_l,w_r)$ defined by
	\begin{equation}
		I^+(\alpha,w_l,w_r)=\frac{1}{\Delta x /2 }\int_{0}^{\Delta x/2} w(\Delta t, x)\di x
		, \quad
		I^-(\alpha,w_l,w_r)=\frac{1}{\Delta x /2 }\int_{-\Delta x/2}^{0} w(\Delta t, x)\di x
	\end{equation}
	If $w$ is seen as a density function, then $I^+(\alpha,w_l,w_r)$ (resp. $I^-(\alpha,w_l,w_r)$) is the value of the density, at time $\Delta t$, of the cell $[0, \Delta x/2]$ (resp. $[-\Delta x/2, 0]$), i.e. the cell of size $\Delta x/2$  whose left (resp. right) border is located at the initial position of the discontinuity (cf. \Cref{fig:sk_cell} for a visual representation).
	
	\begin{figure}
		\centering
		\resizebox{0.6\textwidth}{!}{
			% flatex input: [sk_cell.tex]
			\begin{tikzpicture}[x=\textwidth,y=0.25\textwidth,scale=1, every node/.style={scale=1.2}]

				\draw[thick,dashed,red] (0.2,0.25) -- (0.25,0.25) node[above] {$w_l$};
				\draw[thick,red] (0.25,0.25) -- (0.5,0.25);
				\draw[thick,red](0.5,0.85) -- (0.75,0.85);
				\draw[thick,dashed,red](0.75,0.85) node[above] {$w_r$} -- (0.8,0.85) ;
				\draw[thick,red](0.5,0.85) -- (0.5,0.25);

				\draw[thick,blue, dotted] (0.25,0.25) -- (0.6,0.25);
				\draw[thick,blue, dotted](0.6,0.85) -- (0.75,0.85);
				\draw[thick,blue, dotted](0.6,0.85) -- (0.6,0.25);
				\draw[-latex,blue] (0.5,0.55) -- (0.6,0.55) node[midway, below] {$\alpha \Delta t >0$};
				
				\draw[thick,blue, dotted] (0.25,0.25) -- (0.4,0.25);
				\draw[thick,blue, dotted](0.4,0.85) -- (0.75,0.85);
				\draw[thick,blue, dotted](0.4,0.85) -- (0.4,0.25);
				\draw[-latex,blue] (0.5,0.45) -- (0.4,0.45) node[midway, below] {$\alpha \Delta t <0$};
				
				\draw[dashed] (0.35,1) -- (0.35,-0.1);
				\draw[dashed] (0.5,1) -- (0.5,-0.1);
				\draw[dashed] (0.65,1) -- (0.65,-0.1);

				\draw[-{Latex[length=3mm,width=2mm]}] (0.2,0.15) -- (0.8,0.15) node[anchor=west]{$x$};
				\draw (0.35,0.12) node[anchor=north east]{$-\Delta x/2$} -- (0.35,0.18) ;
				\draw (0.5,0.12) node[anchor=north east] (tot) {$0$} -- (0.5,0.18) ;
				\draw (0.65,0.12) node[anchor=north west]{$\Delta x/2$} -- (0.65,0.18) ;

				\draw[decorate,decoration={brace,mirror}] (0.35,0) -- (0.5,0) node[below,midway,yshift=-5pt] { $I^-(\alpha,w_l,w_r)$};
				
				\draw[decorate,decoration={brace,mirror}] (0.5,-0) -- (0.65,-0) node[below,midway,yshift=-5pt] { $I^+(\alpha,w_l,w_r)$};
				
			\end{tikzpicture}
			
			% flatex input end: [sk_cell.tex]
			
		}
		\caption{Cell averages at the discontinuity. The discontinuity is represented in red for the initial time and in blue after $\Delta t$.}
		\label{fig:sk_cell}
	\end{figure}
	
	Using the expression of the solution to~\eqref{adv} and the condition \eqref{ccfl}, we have
	\begin{equation}
		\begin{aligned}
			I^+(\alpha,w_l,w_r) &=\begin{cases}
				w_r & \text{if } \alpha \le 0 \\
				w_r-2\alpha\frac{\Delta t}{\Delta x} (w_r-w_l) & \text{if } \alpha > 0
			\end{cases}, \\ 
			I^-(\alpha,w_l,w_r) &=\begin{cases}
				w_l - 2\alpha\frac{\Delta t}{\Delta x}(w_r-w_l) & \text{if } \alpha \le 0 \\
				w_l & \text{if } \alpha > 0
			\end{cases}
		\end{aligned}
		\label{it}
	\end{equation}
	Indeed, for instance for $I^+(\alpha,w_l,w_r)$, the condition \eqref{ccfl} gives that $-\Delta x/2 \le -\vert \alpha\vert \Delta t \le 0 \le \vert\alpha\vert \Delta t \le \Delta x/2$. Hence, if $\alpha \ge 0$, we have
	\begin{equation*}
		\begin{aligned}
			I^+(\alpha,w_l,w_r)&=\frac{1}{\Delta x /2 }\left(\int_{0}^{\alpha \Delta t} w(\Delta t, x)\di x+\int_{\alpha \Delta t}^{\Delta x/2} w(\Delta t, x)\di x\right)\\
			&=\frac{1}{\Delta x /2 }\left(\int_{0}^{\alpha \Delta t} w_l \di x+\int_{\alpha \Delta t}^{\Delta x/2} w_r \di x\right) \\
			&=\frac{1}{\Delta x /2 }\left( w_l\alpha \Delta t +\left(\Delta x/2-\alpha \Delta t\right) w_r\right)
			=2\alpha\frac{\Delta t}{\Delta x} w_l+\left(1-2\alpha\frac{\Delta t}{\Delta x}\right) w_r \\
			&=w_r-2\alpha\frac{\Delta t}{\Delta x} (w_r-w_l)
		\end{aligned}
	\end{equation*}
	and since $\alpha \Delta t \ge 0$, for all $x \in [-\Delta x/2, 0[$ we have $x < \alpha\Delta t$, which gives
	\begin{equation*}
		\begin{aligned}
			I^-(\alpha,w_l,w_r)&=\frac{1}{\Delta x /2 }\int_{-\Delta x/2}^{ 0} w_l\di x= w_l.
		\end{aligned}
	\end{equation*}

\section{Additional proofs}

\subsection{Proof of \Cref{thm:genFVM}}\label{sec:proofThm}

\begin{proof}
	%The proposed semi-discretization has the following properties.
	%\begin{itemize}
	%\item 
	The numerical flow $F$ is consistent since by definition and following~\eqref{eq:decomposed_flux}, we have $F(u,u) = g(u, \rho_{\max} -u) = f(u)$ for any $u\in\Omega$. It is monotone due to the continuity assumption on $f$, consistency and the non-decreasing properties of $g$.

	The TRM preserves non-negativity, since, for any $i\in I$,  the right-hand side of~\eqref{eq:gen_scheme}, evaluated on the boundary $\rho_i=0$, is nonnegative. Indeed, following the definitions of $g, g_{\text{on}}, g_{\text{off}}$, if $\rho_i=0$, then
	\begin{equation*}
		\begin{aligned}
			\dot{\rho}_i=\frac{1}{\Delta x}\big[ 
			g(\rho_{i-1}, \rho_{\max}) + R_i(0,t)\big] \geq 0.
		\end{aligned}
	\end{equation*}
	since for any $t\ge 0$, $R_i(\rho_{\max},t)=S_i(0,t)=0$.
	Similarly, the TRM preserves bounded capacity, because the right-hand side of \eqref{eq:gen_scheme}, evaluated on the capacity bound $\rho_i=\rho_{\max}$,  yields $\dot{\rho}_i\le 0$.

	Finally, the discretization scheme \eqref{eq:gen_scheme} is conservative, because for any $n_l, n_r\in I$ such that $n_l < n_r$, the following is true:
	\begin{equation*}
		\Delta x\sum_{i=n_l}^{n_r} \dot{\rho}_i =  F(\rho_{n_l -1}, \rho_{n_l}) - F(\rho_{n_r}, \rho_{n_r+1})
		+\sum_{i=n_l}^{n_r} R_i(\rho_i, t) - S_i(\rho_i, t),
	\end{equation*}
	meaning that the variation of the number of vehicles on a section of the road (i.e. the left hand side of this equation)  is equal to the difference between the number of vehicles entering the section (at the boundary or through on-ramps) and the number of vehicles leaving the section (at the boundary or through off-ramps).
\end{proof}

%--------------------------------------------------------------------------------------

\subsection{Proof of Proposition \ref{prop:eq_pt}}\label{appen:proof_eq}
% flatex input: [Content/proof_equilibrium.tex]

% \begin{prop}
	% For the system of ODEs  resulting from the mass kinetic TRM with a ring topology,  the only possible equilibrium point $(\rho_1^*, \dots, \rho_P^*)$ is the one satsfisfying
	% \begin{equation*}
		%     \rho_i^* = \bar{\rho}=\frac{1}{P}\sum_{i=1}^P\rho_i(0)
		% \end{equation*}
	% \end{prop}

\begin{proof}
	
	Since $(\rho_1^*, \dots, \rho_P^*)$ is an equilibrium point, it must satisfy
	\begin{equation*}
		F(\rho_{i-1}^*, \rho_i^*)=F(\rho_{i}^*, \rho_{i+1}^*), \quad i\in \lbrace 1, \dots, P\rbrace
	\end{equation*}
	where we write in particular $\rho_0^*=\rho_P^*$ and $\rho_{P+1}^*=\rho_1^*$.
	This gives in turn,
	\begin{equation*}
		\rho_P^* (\rho_{\max} - \rho_{1}^*) = \rho_1^* (\rho_{\max} - \rho_2^*) = \dots = \rho_{P-1}^* (\rho_{\max} - \rho_P^*).
	\end{equation*}
	Without loss of generality, assume that $\rho_1^*=\max_{i\in I} \rho_i^*$. We can then write,
	\begin{equation*}
		\rho_1^*(\rho_{\max}-\rho_2^*)-\rho_1^*(\rho_{\max}-\rho_1^*)
		=\rho_P^*(\rho_{\max}-\rho_1^*)-\rho_1^*(\rho_{\max}-\rho_1^*)
	\end{equation*}
	which gives
	\begin{equation*}
		\rho_1^*(\rho_1^*-\rho_2^*)=(\rho_{\max}-\rho_1^*)(\rho_P^*-\rho_1^*)
	\end{equation*}
	and therefore
	\begin{equation*}
		\underbrace{\rho_1^*}_{\ge 0}\underbrace{(\rho_1^*-\rho_2^*)}_{\ge 0}+\underbrace{(\rho_{\max}-\rho_1^*)}_{\ge 0}\underbrace{(\rho_1^*-\rho_P^*)}_{\ge 0}=0
	\end{equation*}
	We have three cases.
	\begin{itemize}
		\item If $\rho_1^*=0$, then since $\rho_1^*=\max_{i\in I} \rho_i^*$, we have $\rho_1^*=\rho_2^*=\dots=\rho_P^*=0$.
		\item If $\rho_1^*=\rho_{\max}$, then since $\rho_1^* (\rho_{\max} - \rho_2^*)=\rho_P^* (\rho_{\max} - \rho_{1}^*) =0$, we have $\rho_2^*=\rho_{\max}$. And  since $\rho_2^* (\rho_{\max} - \rho_3^*)=\rho_1^* (\rho_{\max} - \rho_{2}^*) =0$, we have $\rho_3^*=\rho_{\max}$. By continuing this process, we then prove that $\rho_1^*=\rho_2^*=\dots=\rho_P^*=\rho_{\max}$.
		\item If $\rho_1^* \neq 0$ and $\rho_1^* \neq \rho_{\max}$: then we must have $\rho_P^*=\rho_1^*=\rho_2^*=\max_{i\in I} \rho_i^*$. Hence we can apply the same result to prove that $\rho_1^*=\rho_2^*=\rho_3^*=\max_{i\in I} \rho_i^*$, and continue to do iteratively until we prove that $\rho_{P-2}^*=\rho_{P-1}^*=\rho_P^*=\max_{i\in I} \rho_i^*$. Hence we have $\rho_1^*=\rho_2^*=\dots=\rho_P^*$.
	\end{itemize}
	In all three cases we proved that $\rho_1^*=\rho_2^*=\dots=\rho_P^*$, which is clearly equal to $(1/P)\sum_{i\in I}\rho_i^*$. 
	
	Note then that the quantity
	\begin{equation*}
		\bar{\rho}(t)=\frac{1}{P}\sum_{i=1}^P\rho_i(t)
	\end{equation*}
	satisfies
	\begin{equation*}
		\dot{\bar{\rho}}(t)=(F(\rho_P,\rho_1)-F(\rho_P,\rho_1))/P=0
	\end{equation*}
	since the TRM is conservative. Hence, $\bar{\rho}(t)$ is conserved through time, meaning that $\bar{\rho}(t)=(1/P)\sum_{i\in I}\rho_i^*=\bar{\rho}(0)=\bar{\rho}$. We can then conclude that $\rho_1^*=\rho_2^*=\dots=\rho_P^*=\bar{\rho}$.
	
\end{proof}

% flatex input end: [Content/proof_equilibrium.tex]

\subsection{Proof of Proposition \ref{prop:modeq_shock}}\label{sec:proof_meq}
% flatex input: [Content/proof_meq.tex]

\begin{proof}
	Let $\tilde\rho : \R_+\times\R \rightarrow \R$ be the solution of the following PDE
	\begin{equation}
		\partial_t\tilde{\rho} + \omega \partial_x\big(\tilde{\rho}(\rho_{\max}-\tilde{\rho})\big) + \frac{\Delta t}{2}\partial_{tt}\tilde{\rho}-\frac{\Delta x}{2}\omega\rho_{\max}\partial_{xx}\tilde{\rho}=0, 
	\end{equation}
	and whose initial condition $\rho(0, \cdot)$ satisfies
	\begin{equation}
		\lim\limits_{x\rightarrow -\infty} \rho(0, x)=\rho_l, \quad
		\lim\limits_{x\rightarrow +\infty} \rho(0, x)=\rho_r
		\label{eq:modproof2}
	\end{equation}
	where $\rho_l < \rho_r$.
	
	Let us consider traveling wave solutions to this PDE, meaning that we assume that $\tilde{\rho}$ can be written as 
	\begin{equation}
		\tilde{\rho}(t,x)=f(x-ct), \quad t\ge 0, \quad x\in \R, 
		\label{eq:modproof1}
	\end{equation}
	for some (at least twice-differentiable) function $f$ and $c\in\R$. In this case, note that \eqref{eq:modproof2} yields the following requirement for $f$:
	\begin{equation}
		\lim\limits_{x\rightarrow -\infty} f(x)=\rho_l, \quad
		\lim\limits_{x\rightarrow +\infty} f(x)=\rho_r.
		\label{eq:modproof2b}
	\end{equation}
	
	Then, the PDE becomes:
	\begin{equation}
		\big(\frac{\Delta t}{2}c^2-\omega\rho_{\max}\frac{\Delta x}{2}\big)f''(z) + \big(\omega\rho_{\max}-c-2\omega f(z)\big)f'(z)=0, \label{eq:modproof3}
	\end{equation}
	where $z=x-ct \in\R$.
	
	Let $A, B_1, B_2$ be the constants defined by
	\begin{equation*}
		A=\frac{\Delta t}{2}c^2-\omega\rho_{\max}\frac{\Delta x}{2}, \quad B_1=2\omega, \quad B_2=c-\omega\rho_{\max}.
	\end{equation*}
	Note that if $A=0$ then we  get from \eqref{eq:modproof3} that $f$ is equal to a constant for any $z\in\R$ (since it is twice-differentiable), and therefore a contradiction with \eqref{eq:modproof2b}. Hence, can assume that $A\neq 0$.
	Consider then the change of variables defined by
	\begin{equation}
		f(z)=-\frac{2A}{B_1} \frac{\psi'(z)}{\psi(z)}, \quad z\in\R,
		\label{eq:modproof6}
	\end{equation}
	where (without loss of generality) $\psi(z)>0$ for any $z\in\R$. Then, \eqref{eq:modproof3} becomes,
	\begin{equation*}
		-\frac{2A^2}{B_1}\frac{\psi'''\psi-\psi'\psi''}{\psi^2}+\frac{2B_2A}{B_1}\frac{\psi''\psi-(\psi')^2}{\psi^2}=0
	\end{equation*}
	which gives
	\begin{equation*}
		\frac{B_2\psi''-A\psi'''}{\psi}-\frac{B_2(\psi')^2-A\psi'\psi''}{\psi^2}=0
	\end{equation*}	
	and therefore
	\begin{equation*}
		\frac{d}{dz}\bigg[ \frac{B_2\psi'(z)-A\psi''(z)}{\psi(z)} \bigg]=0.
	\end{equation*}
	Hence, $\psi$ must satisfy an ODE of the form:
	\begin{equation}
		\psi''(z)= \frac{B_2}{A} \psi'(z) + K\psi(z), \quad z\in\R.
		\label{eq:modproof6b}
	\end{equation}
	for some constant $K\in\R$. 
	
	The solutions of this ODE depend on the sign of the quantity $\lambda$ defined by
	\begin{equation}
		\lambda = \bigg(-\frac{B_2}{A}\bigg)^2+4K.
	\end{equation}
	Since $\psi>0$ and therefore is never zero, we must have $\lambda >0$ and so, $4K>-(B_2/A)^2$. In that case, the solutions of \eqref{eq:modproof6b} can be written as
	\begin{equation}
		\psi(z)=D_1 e^{r_1z/2}+D_2 e^{r_2z/2}, \quad z\in\R, 
	\end{equation}
	for some constants $D_1, D_2\in\R$ and where $r_1, r_2$ are defined as
	\begin{equation*}
		r_1=\frac{B_2}{A}+\sqrt{\lambda}, \quad r_2=\frac{B_2}{A}-\sqrt{\lambda}.
	\end{equation*}
	%Note in particular that then, $r_1r_2 <0$. Hence, $r_1$ and $r_2$ are non-zero and have different signs. 
	Once again, since $\psi >0$,  $D_1$ and $D_2$ cannot be identically equal to $0$, and must be non-negative (otherwise one can find some $z\in\R$ such that $\psi(z)=0$).
	
	In turn,  $f$ may be written as
	\begin{equation*}
		f(z)=-\frac{A}{B_1}
		\frac{D_1r_1 e^{r_1z/2}+D_2r_2 e^{r_2z/2}}{D_1e^{r_1z/2}+D_2e^{r_2z/2}}, \quad z\in\R.
	\end{equation*}
	Nothing that $r_1>r_2$, we rewrite $f$ as
	\begin{equation*}
		f(z)=-\frac{A}{B_1}
		\frac{D_1r_1 e^{(r_1-r_2)z/2}+D_2r_2}{D_1e^{(r_1-r_2)z/2}+D_2}, \quad z\in\R.
	\end{equation*}

	Taking then the limit at $\pm\infty$ of this expression, we see that in order to get two different limits (namely $\rho_l$ and $\rho_r$) we must have that both $D_1$ and $D_2$ are non-zero. We then have
	\begin{equation*}
		\rho_r=-\frac{A}{B_1}r_1, \quad \rho_l=-\frac{A}{B_1}r_2.
	\end{equation*}
	This gives, on the one hand,
	\begin{equation*}
		r_1-r_2=-\frac{B_1}{A}(\rho_r-\rho_l)=-\frac{2\omega}{A}(\rho_r-\rho_l).
	\end{equation*}
	Hence, since $r_1-r_2>0$ and $\rho_r -\rho_l>0$ we have that $A<0$ and we can write
	\begin{equation*}
		f(z)=
		\frac{D_1\rho_r e^{-\omega(\rho_r-\rho_l) z/A}+D_2\rho_l}{D_1e^{-\omega (\rho_r-\rho_l)z/A}+D_2}, \quad z\in\R.
	\end{equation*}
	On the other hand,	
	\begin{equation*}
		\rho_l+\rho_r=-\frac{A}{B_1}(r_1+r_2)=-2\frac{B_2}{B_1}=\frac{\omega \rho_{\max}-c}{\omega}.
	\end{equation*}
	Hence, we have
	\begin{equation}
		c=\omega(\rho_{\max}-\rho_l-\rho_r).
	\end{equation}
	This gives on the one hand,
	\begin{equation*}
		B_2=-\omega(\rho_l+\rho_r)<0,
	\end{equation*}
	and on the other hand,
	\begin{equation*}
		\begin{aligned}
			A&=\frac{\Delta t}{2}(\omega(\rho_{\max}-\rho_l-\rho_r))^2-\omega\rho_{\max}\frac{\Delta x}{2} \\
			&=\frac{\omega\rho_{\max}}{2}\Delta x \bigg(\omega\rho_{\max}\frac{\Delta t}{\Delta x}\bigg(1-\frac{\rho_l+\rho_r}{\rho_{\max}}\bigg)^2-1\bigg)
		\end{aligned}
	\end{equation*}
	Note then that if we take $\delta >0$ as
	\begin{equation*}
		\delta = \omega\rho_{\max}\frac{\Delta t}{\Delta x},
	\end{equation*}
	according to the CFL condition, we have that $\delta < 1/2$, and since $\rho_l, \rho_r \in [0, \rho_{\max}]$, we retrieve that
	\begin{equation*}
		A
		=-\frac{\omega\rho_{\max}}{2}\Delta x \bigg(1-\delta\bigg(1-\frac{\rho_l+\rho_r}{\rho_{\max}}\bigg)^2\bigg) < 0.
	\end{equation*}

	Consequently, we can write
	\begin{equation*}
		f(z)=
		\frac{D_1\rho_r e^{z/(2\sigma)}+D_2\rho_l e^{-z/(2\sigma)}}{D_1 e^{z/(2\sigma)}+D_2 e^{-z/(2\sigma)}}
		=\rho_l + \frac{\rho_r-\rho_l}{1+De^{-z/\sigma}}
		,
	\end{equation*}
	where $z\in\R$ and $D=D_2/D_1 >0$ and
	\begin{equation}
		\sigma=-\frac{A}{\omega(\rho_r-\rho_l)}=\frac{1}{2}\frac{\rho_{\max}}{\rho_r-\rho_l} \bigg(1-\delta\bigg(1-\frac{\rho_l+\rho_r}{\rho_{\max}}\bigg)^2\bigg) \Delta x.
	\end{equation}
	Finally, note that since $D>0$ we can take $x_0=\log(D)\in\R$ to retrieve 
	\begin{equation}
		f(z)
		=\rho_l + \frac{\rho_r-\rho_l}{1+e^{-(z-x_0)/\sigma}}
		, \quad z\in\R.  
	\end{equation}		
\end{proof}

% flatex input end: [Content/proof_meq.tex]

%--------------------------------------------------------------------------------------

\end{appendices}
%%---------------------------------------------------------------------------

\end{document}